\documentclass[11pt]{article}
\usepackage{amsmath,amssymb,latexsym}
\usepackage{theorem}
\usepackage[final]{graphicx}
\usepackage{overpic}
\newtheorem{theorem}{Theorem}

\newtheorem{corollary}{Corollary}
\newtheorem{proposition}{Proposition}
{\theorembodyfont{\rmfamily} \newtheorem{definition}{Definition}

\newtheorem{remark}{Remark}} {\theorembodyfont{\slshape} }

\newcommand{\res}{\mathop{\rm res}}

\newcommand{\mult}{\mathop{\rm mult}}
\newcommand{\field}[1]{\mathbb{#1}}
\newcommand{\D}{\field{D}}
\newcommand{\R}{\field{R}}
\newcommand{\Z}{\field{Z}}
\newcommand{\N}{\field{N}}
\newcommand{\C}{\field{C}}
\newcommand{\Q}{\field{Q}}

\newcommand{\T}{\field{T}}
\newcommand{\CC}{{\mathcal C}}
\newcommand{\DD}{{\mathcal D}}
\newcommand{\EE}{{\mathcal E}}
\newcommand{\II}{{\mathcal I}}

\newcommand{\FF}{{\mathcal F}}
\newcommand{\PP}{{\mathcal P}}

\renewcommand{\Re}{\mathop{\rm Re}}
\renewcommand{\Im}{\mathop{\rm Im}}
\newcommand{\dist}{\mathop{\rm dist}}
\newcommand{\isdef}{\stackrel{\text{\tiny def}}{=}}

\DeclareRobustCommand{\qed}{%
\ifmmode 
\else \leavevmode\unskip\penalty9999 \hbox{}\nobreak\hfill \fi
\quad\hbox{\qedsymbol}}
\newcommand{\openbox}{\leavevmode
\hbox to.77778em{%
\hfil\vrule
\vbox to.675em{\hrule width.6em\vfil\hrule}%
\vrule\hfil}}
\newcommand{\qedsymbol}{\openbox}
\newcommand{\proofname}{Proof}
\newenvironment{proof}[1][\proofname]{\par
\normalfont \trivlist \item[\hskip\labelsep   \itshape #1. ]
\ignorespaces
}{%
\qed\endtrivlist } 

\def\Xint#1{\mathchoice
{\XXint\displaystyle\textstyle{#1}}%
{\XXint\textstyle\scriptstyle{#1}}%
{\XXint\scriptstyle\scriptscriptstyle{#1}}%
{\XXint\scriptscriptstyle\scriptscriptstyle{#1}}%
\!\int}
\def\XXint#1#2#3{{\setbox0=\hbox{$#1{#2#3}{\int}$}
\vcenter{\hbox{$#2#3$}}\kern-.5\wd0}}

\def\dashint{\Xint-}

\title{Szeg\H{o} orthogonal polynomials with respect to an analytic weight: canonical representation
and strong asymptotics.%
}%
\author{A.\ Mart\'{\i}nez-Finkelshtein\footnote{Corresponding author. Email: \texttt{andrei@ual.es}}\\ Universidad de Almer\'{\i}a, Spain
\and K. T.-R. McLaughlin\\ University of Arizona, USA \and
E.\ B.\ Saff \\Vanderbilt University, USA }%
\begin{document}
\maketitle

\begin{abstract}
We provide a representation in terms of certain canonical functions
for a sequence of polynomials orthogonal with respect to a weight
that is strictly positive and analytic on the unit circle. These
formulas yield a complete asymptotic expansion for these
polynomials, valid uniformly in the whole complex plane. As a
consequence, we obtain some results about the distribution of zeros
of these polynomials. The main technique is the steepest descent
analysis of Deift and Zhou, based on the matrix Riemann-Hilbert
characterization proposed by Fokas, Its and Kitaev.
\end{abstract}

\section{Introduction and background} \label{section:the weight}

For $r>0$, let $\D_r$ be the open disc in $\C$ with center at
$z=0$ and radius $r$, $\{z\in \C:\, |z|<r\}$, and let $\T_r$ be
its boundary $ \{z\in \C:\, |z|=r \}$. An integrable non-negative
function $w$ defined on $\T_1$ is called a \emph{weight} if
$$
\int_{\T_1} w(z)\, |dz| >0\,.
$$
For each weight $w$ there exists a unique sequence of polynomials
$\varphi_n$ (called \emph{Szeg\H{o} polynomials}), orthonormal
with respect to $w$, satisfying $\varphi_n (z)=\kappa_n
z^n+\text{lower degree terms}$, $\kappa_n >0$, and
\begin{equation}\label{orthogonalityConditions}
\oint_{\T_1} \varphi_n(z)  \overline{\varphi_m(z)}\, w(z)
|dz|=\delta_{mn}\,.
\end{equation}
We denote by $\Phi_n(z) \isdef  \varphi_n(z)/\kappa _n$ the
corresponding monic orthogonal polynomials. It is well known that
they satisfy the Szeg\H{o} recurrence
$$
\Phi_{n+1}(z)=z \Phi_n(z)-\overline{\alpha _n}\, \Phi_n^*(z)\,,
\quad \Phi_0(z)\equiv 1\,,
$$
where we use the standard notation $\Phi_n^*(z)\isdef z^n
\overline{\Phi_n(1/\overline{z})}$. The parameters $\alpha
_n=-\overline{\Phi_{n+1}(0)}$ are called \emph{Verblunsky
coefficients} (also \emph{reflection coefficients} or \emph{Schur
parameters}) and satisfy $\alpha _n \in \D_1$ for $n=0, 1, 2, \dots$
(see \cite{Simon04a} for details).

The asymptotics of the sequence $\{\varphi_n (z)\}$ when $n\to
\infty$ is well known for $|z|>1$  (the so called \emph{outer} or
\emph{exterior asymptotics}) under rather general assumptions on
the weight. However, the behavior of these polynomials in $\D_1$
is, in general, much more complicated and less studied. First, it
is known that all the zeros of each $\Phi_n$ lie in $\D_1$ (see
e.g.\ \cite[p.\ 50--55]{Grenander/Szego:1984}, \cite[Ch.\ XII,
theorem 12.1.1]{szego:1975}, \cite{Geronimus61} or
\cite{Simon04a}). Moreover, they can be everywhere dense in the
unit disc (the so-called Turan's conjecture), as shown in
\cite{Alfaro88}. A recent result of Simon and Totik
\cite{Simon/Totik:04} proves even the existence of a ``universal''
measure of orthogonality on $\T_1$ for which the zero
distributions of the corresponding $\{\Phi_n (z)\}$ approximate
(in a weak sense) all the possible measures on the closed disc
$\overline{\D}_1$.

However, under stronger assumptions on the weight $w$ we can say
much more. In what follows we will consider only the strictly
positive weights $w$ on the unit circle that can be extended to an
open set containing $\T_1$ as a nowhere vanishing holomorphic
function. We will say briefly that such $w$ is a \emph{strictly
positive analytic weight} on $\T_1$. All the forthcoming
statements are valid under this assumption. In order to shed light
on this definition we point out that  $w$ is a non-vanishing
holomorphic function in an annulus $r <|z|<1/r $ (with $0<r<1$)
and strictly positive on the unit circle $\T_1$, if and only if,
for any nonzero $\alpha\in \R$ the following representation is
valid on $\T_1$:
\begin{equation}\label{on the circle}
w(z)=|W(z)|^\alpha\,, \quad z \in \T_1\,,
\end{equation}
where $W$ is holomorphic and nowhere vanishing in the annulus. In
particular, for $\alpha =2$, the following representation is
useful:
\begin{equation}\label{representation weight}
w(z)=W(z) \overline{W(z)}=W(z) \overline{W(1/\overline{z})}\,,
\quad z\in \T_1\,.
\end{equation}
Obviously, here $W$ is not uniquely defined. Normalized
representations of this form arise from the so-called
\emph{Szeg\H{o} function} of $w$. Indeed, since $w$ is strictly
positive on $\T_1$, its index (winding number about the origin) is
$0$, and we may define the Szeg\H{o} function (see e.g.\
\cite[Ch.\ X, \S 10.2]{szego:1975}):
\begin{equation}\label{standardSzego}
D (z) \isdef \exp\left( \frac{1}{4\pi }\,\int_0^{2\pi} \log w(e^{i
\theta}) \, \frac{e^{i \theta}+z}{e^{i \theta}-z}\, d\theta
\right)\,.
\end{equation}
This function is piecewise analytic and non-vanishing, defined for
$|z|\neq 1$, and we will denote by $D_{\rm i}$ and $D_{\rm e}$ its
values for $|z|<1$ and $|z|>1$, respectively. It is easy to verify
that
\begin{equation}\label{symmetry}
\overline{D_{\rm i}\left(\frac{1}{\overline{z}}
\right)}=\frac{1}{D_{\rm e}(z)}\,, \quad |z|>1\,.
\end{equation}
Furthermore, if $w$ satisfies the assumptions above, both $D_{\rm
i}$ and $D_{\rm e}$ admit a holomorphic extension across $\T_1$,
and maintaining the same notation for these analytic continuations
we have
\begin{equation}\label{boundary_value_Szego}
\frac{D_{\rm i}(z)}{D_{\rm e}(z)}=w(z)\,.
\end{equation}
This formula gives an analytic extension of $w$ from $\T_1$ to an
open domain containing the unit circle and, up to a normalization
constant, provides the ``standard'' representation
\eqref{representation weight}. Observe that by \eqref{symmetry},
$w(z)=|D_{\rm e}(z)|^{-2}$ for $z\in \T_1$.

Let us define
\begin{equation}\label{Def_R_bis}
\rho \isdef  \inf\{0<r<1:\, D_{\rm e}(z) \text{ is holomorphic in
} |z|>r \}\,.
\end{equation}
The class of weights that we consider in this article is
characterized by the fact that $\rho <1$. Due to the symmetry
\eqref{symmetry}, $1/\rho $ is the radius of convergence of the
Taylor series for $1/D_{\rm i}$ about $z=0$ (hence, the constant
$1/\rho $ coincides with $R_{D^{-1}}$ in the notation of \cite[\S
7.1]{Simon04a}). Taking into account \eqref{boundary_value_Szego},
we may equivalently define $\rho$ as
\begin{equation}\label{Def_R}
\rho=\inf\{0<r<1:\, 1/w(z) \text{ is holomorphic in } r<|z|<1/r
\}\,.
\end{equation}
Again by \eqref{symmetry}, both circles $\T_\rho $ and
$\T_{1/\rho} $ contain singularities  of $1/w$.

Important results concerning the class of strictly positive
analytic weights on $\T_1$ were obtained in \cite{Nevai/Totik:89}.
It was shown there (see also \cite[\S 7.1]{Simon04a}) that
\begin{equation}\label{resultNevaiTotik}
\rho=\varlimsup_{n\to \infty} |\alpha _n|^{1/n}\,.
\end{equation}
In particular, strictly positive analytic weights on $\T_1$ are
characterized by an exponential decay of the corresponding
Verblunsky coefficients. Furthermore, the Szeg\H{o} asymptotic
formula
\begin{equation}\label{ExteriorSzego}
\lim_{n\to \infty}\varphi_n^* (z) =\frac{1}{D_{\rm i}(z)}\,, \quad
z\in \D_1\,,
\end{equation}
can be continued analytically through the unit circle $\T_1$ and is
valid locally uniformly in $\D_{1/\rho }$. It shows that the number
of zeros of $\{\varphi _n\}$ outside of $\overline{\D}_\rho $
remains uniformly bounded and these zeros are attracted by the zeros
of $D_{\rm e}$ in $\rho <|z|<1$ (\emph{Nevai-Totik points} in the
terminology of B.\ Simon \cite{Simon04b}). Mhaskar and Saff
\cite{Mhaskar90}, using potential theory arguments, gave a relevant
complement to these results about zeros showing that for any
subsequence $\{ n_k \} \subset \N$ satisfying
$$
\rho=\lim |\alpha _{n_k}|^{1/n_k}\,,
$$
the zeros of $\{\varphi_{n_k+1}\}$ distribute asymptotically
uniformly in the weak-star sense on $\T_\rho $ (see also
\cite{Pakula87}).

Further results on the asymptotics or zero behavior inside
$\T_\rho $ have been obtained so far only for some particular
subclasses of strictly positive analytic weights on $\T_1$.
Besides the classical explicit formula for the Bernstein-Szeg\H{o}
polynomials (when $w(z)=|1/q(z)|^2$, $q$ a fixed polynomial with
zeros in $\C\setminus \T_1$), probably the first important result
in this direction was obtained by Szabados \cite{Szabados79}, who
considered the polynomial case $w(z)=|p(z)|^2$ (where again $p$ is
a fixed polynomial with zeros in $\C\setminus \T_1$). Szabados
gave a formula for the orthogonal polynomials with respect to
$|p(z)|^2$ which, although containing unknown parameters, is
sufficient for the asymptotic analysis inside $\T_\rho $.

By extending Darboux's formula to Szeg\H{o} orthogonal polynomials,
Ismail and Ruedemann \cite{Ismail92} (see also \cite{Godoy91})
obtain matrix expressions for the $\varphi _n$'s for the case of
rational weights $w(z)=|p(z)/q(z)|^2$, where the size of the matrix
depends on the degrees of $p$ and $q$. They applied their result to
the asymptotic analysis in the case when $\T_\rho $ contains a
single zero of $w$. The same rational case is also the object of
study in a less known paper \cite{Pakula87}, where the author
considers the weak-star zero distribution for the $\varphi _n$'s.

Another class of weights for which $\rho <1$ corresponds to the
Rogers-Szeg\H{o} polynomials, studied for instance, in
\cite{Mazel/Geronimo/Hayes:90}, where even the interlacing
properties of their zeros have been established (see also
\cite{Lubinsky/Saff87} and \cite[Example 1.6.5]{Simon04a}). The
reader interested in many more explicit examples of orthogonal
polynomials on $\T_1$ is referred to \cite[Ch.\ I and
VIII]{Simon04a}.

The main goal of this paper is to provide a complete asymptotic
expansion for the orthogonal polynomials, valid uniformly in the
whole complex plane, under the assumption that $w$ is a strictly
positive analytic weight on $\T_1$. As a consequence, we obtain
results about the distribution of zeros of these orthogonal
polynomials, under a variety of different assumptions about the
weight function.  These results recover, in particular cases, some
of the results mentioned above. The main technique is the steepest
descent analysis of Deift and Zhou (see e.g.\ \cite{MR2000g:47048}),
based on the matrix Riemann-Hilbert characterization proposed by
Fokas, Its and Kitaev \cite{Fokas92}. This approach more generally
allows one to handle analytic complex valued weights $w$ that are
nonvanishing on $\T_1$ provided the corresponding orthonormal
polynomials $\varphi_n (z)$ exist with $\kappa_n \neq 0$ for each
$n$. Furthermore, the asymptotic analysis is equally possible (and
technically not much more involved) for a sequence of varying
weights (when $w$ depends on the degree $n$ of the polynomial). This
was carried out previously in \cite{MR2000e:05006} for a specific
weight relevant to the asymptotics of the distribution of the length
of the longest increasing sequences in random permutations. A basic
discussion of orthogonal polynomials on the unit circle, in the
context of an example of integrable operators, and their connection
with the Riemann-Hilbert problems is contained in \cite{deift99}.

Until recently the analyticity of $w$ has been a necessary
condition for the applicability of the Riemann-Hilbert asymptotic
analysis. However, McLaughlin and Miller
\cite{McLaughlin/Miller:2004} have extended the approach via the
solution of a d-bar problem. Their method allows one to handle
piecewise smooth and strictly positive weights $w$ on $\T_1$; for
this generality they pay the price of less detailed asymptotics.
Also their results are qualitatively different: in the settings of
\cite{McLaughlin/Miller:2004} essentially all the zeros approach
$\T_1$ with a polynomial rate depending on the smoothness of $w$,
and those that might stay inside are artifacts of the jumps of the
weight.

Our asymptotic formulas (Theorem \ref{thm:nice case}) are stated
in terms of a sequence of iterated Cauchy transforms of the
function
\begin{equation}\label{def_F}
    \FF(z) \isdef D_{\rm i}(z)D_{\rm e}(z)\,,
\end{equation}
where $D_{\rm i}$ and $D_{\rm e}$ are respectively the inner and
the outer Szeg\H{o} functions defined by \eqref{standardSzego}. By
\eqref{Def_R_bis}, $\FF$ is holomorphic in $\rho <|z|<1 $.
Furthermore, by \eqref{symmetry},
\begin{equation}\label{symmetryOfF}
 \overline{\FF\left(\frac{1}{\overline{z}}
\right)}=\frac{1}{\FF(z)}\,, \quad \text{for } |z|\neq 1\,, \qquad
 |\FF(z)|=1 \text{ on } \T_1\,,
\end{equation}
and $1/\FF$ is holomorphic in $1 <|z|<1/\rho $.

The function $1/\FF$ appears as a scattering matrix (or function)
in \cite[Formula (V.5)]{Geronimo/Case:79}. Geronimo and Case also
give an asymptotic expression for $\varphi_n$ on $\T_1$, where the
Cauchy transform of the scattering functions appears implicitly
(cf.\ formulas (III.7) and (V.8) in \cite{Geronimo/Case:79}). One
of our goals here is to make precise the relevant role of the
sequences of iterated Cauchy transforms related to $\FF$ in the
asymptotics of the orthogonal polynomials.

The structure of the paper is as follows. In the next section we
state the main results giving a representation of the polynomials
$\Phi_n$, the leading coefficients $\kappa_{n}$ and the Verblunsky
coefficients $\alpha_{n}$ in terms of the series of some canonical
functions. These formulas have an asymptotic nature that we exploit
in Section \ref{sec:asymptotics}, where we present a number of
results concerning the asymptotic behavior of the Verblunsky
coefficients, leading coefficients, and the ``outer'' asymptotic
behavior of the orthogonal polynomials.  We also present in Section
\ref{sec:asymptotics} general results describing the asymptotic
behavior of the zeros, for a rather wide class of analytic weights.
The detailed assumption is that the exterior function $D_{\rm e}$
can be extended to the exterior of a domain $\{ |z| > \rho' \}$ with
finitely many poles.  The proofs of these results are deferred to
Section \ref{sec:RH_analysis}. Finally, in Section
\ref{sec:examples} we discuss some particular cases, including the
case of a weight for which $D_{e}$ possesses an isolated essential
singularity within the unit disk. The examples in Section
\ref{sec:examples} are all illustrated also with numerical
experiments. We should point out that the zeros of $\Phi_n$'s have
been computed as the eigenvalues of the truncated matrix
corresponding to the CMV representation of the multiplication
operator (see \cite{Cantero03} or \cite[Ch.\ IV]{Simon04a} for
details), using a symbolic algebra software with extended precision.

\section{Canonical representation of the orthogonal polynomials} \label{section:statement}

Occasionally it is more natural to write our formulas in terms of
a modified Szeg\H{o} function.  Since $w$ is of index $0$, we may
define the calligraphic functions
\begin{align}
\label{iterior} \DD_{\rm i}(z)&\isdef \exp\left( \frac{1}{2\pi
i}\,\oint_{\T_1} \frac{\log w(t)\, dt}{t-z}\right),
& \text{if } |z|<1\,, \\
\label{exterior} \DD_{\rm e}(z)&\isdef \exp\left( -\frac{1}{2\pi
i}\,\oint_{\T_1} \frac{\log w(t)\, dt}{t-z}\right), & \text{if }
|z|>1\,.
\end{align}
Both $\DD_{\rm i}$ and $\DD_{\rm e}$ are holomorphic and
non-vanishing in the respective regions where they are defined and
satisfy
\begin{equation}\label{symmetryD}
\tau\,\overline{ \DD_{\rm i}\left(\frac{1}{\overline{z}}
\right)}=\frac{\DD_{\rm e}(z)}{\tau}\,, \qquad |z|>1\,,
\end{equation}
where $\tau $ is the Szeg\H{o} constant (mean value) given by
 \begin{equation}\label{tau}
\tau \isdef \frac{1}{D_{\rm i}(0)}=\exp\left( - \frac{1}{4\pi
}\,\int_0^{2\pi} \log w(e^{i \theta}) \, d\theta
\right)=\exp\left( - \frac{1}{4\pi i }\,\int_{{\T_1}} \log w(t) \,
\frac{dt }{t} \right)>0\,.
\end{equation}
Furthermore (compare \eqref{boundary_value_Szego}),
\begin{equation}\label{boundary_value}
\DD_{\rm i}(z) \DD_{\rm e}(z)=w(z) \,.
\end{equation}
This equation plus the normalization $\DD_{\rm e}(\infty)=1$
determines uniquely the nowhere vanishing functions
\eqref{iterior}--\eqref{exterior}. Finally, we have the following
straightforward identities for $|z|\neq 1$:
\begin{align*}
D_{\rm i}(z) &=  \tau\, \DD_{\rm i}(z),  &\text{if } |z|<1\,, \\
 D_{\rm e}(z) &=  \tau\,   \DD_{\rm e}(z)^{-1}, &\text{if }
 |z|>1\,.
\end{align*}

In what follows, all the circles $\T_\alpha $, $\alpha >0$, are
oriented counterclockwise. For a fixed value of $\alpha$ (which will
be clear from the context), we set
$$
f_+ (z) \isdef \lim_{\stackrel{u\to z}{|u|<\alpha }}f(t)\,, \qquad
f_- (z) \isdef \lim_{\stackrel{u\to z}{|u|>\alpha }}f(t)\,,
$$
for any analytic function $f$ for which these limits exist. We
will also talk about the ``$+$'' and the ``$-$'' side of
$\T_\alpha $ referring to its inner and outer boundary points,
respectively.

Choose an arbitrary $r$, $\rho <r<1$, that we fix for what
follows. We introduce the following operators acting on the space
of holomorphic functions in $\C \setminus ({\T_r} \cup {\T_{1/r}}
)$ with continuous boundary values:
\begin{align} \label{MPlus}
\mathcal M_n^{\rm i}(f)(z)& \isdef -\frac{1}{2\pi i}\,
\oint_{{\T_r}} f_-(t)\, \frac{t^n \DD_{\rm i}^2(t)}{ w(t)(t-z)} \,
dt=-\frac{1}{2\pi i\, \tau^2}\, \oint_{{\T_r}} f_-(t)\,
 \frac{\FF(t) \, t^n }{ t-z} \, dt\,, \\ \label{MMinus}
 \mathcal M_n^{\rm e}(f)(z) & \isdef \frac{1}{2\pi i}\,
\oint_{{\T_{1/r}}} f_-(t)\, \frac{\DD_{\rm e}^2(t)}{ t^n
w(t)(t-z)} \, dt = \frac{\tau^2 }{2\pi i}\, \oint_{{\T_{1/r}}}
f_-(t)\, \frac{ 1 }{\FF(t)\, t^n\, (t-z)} \, dt\,,
\end{align}
where $f_-$ denotes the exterior boundary values of the function
$f$ on the respective circles.

Let us denote by
\begin{equation}\label{defLambda}
\Lambda \isdef \|\FF\|_{\T_r}=\max_{z\in \T_r} |\FF(z)|\,;
\end{equation}
taking into account \eqref{symmetryOfF}, also
$$
\Lambda =\left\|\frac{1}{\FF}\right\|_{\T_{1/r}}=\max_{z\in
\T_{1/r}} \frac{1}{|\FF(z)|}\,.
$$
Then straightforward bounds show that
\begin{equation}\label{bounds_operators}
\begin{split}
\left| \mathcal M_n^{\rm i}(f)(z) \right| & \leq \frac{\Lambda
r}{\tau ^2}\, r^n\, \frac{ \|f_-\|_{\T_r}}{\left||z|-r \right|}\,,
\qquad
z\notin \T_r\,, \\
 \left| \mathcal M_n^{\rm e}(f)(z) \right| & \leq
\frac{\Lambda \tau ^2}{r}\, r^n\, \frac{
\|f_-\|_{\T_{1/r}}}{\left||z|-1/r \right|}\,, \qquad z\notin
\T_{1/r}\,.
\end{split}
\end{equation}

\begin{remark}
It is obvious that for $|z|<\rho $ in \eqref{MPlus} and for
$|z|>1/\rho $ in \eqref{MMinus} we may take $\T_1$ as the contour
of integration in these formulas. However, our choice makes
inequalities \eqref{bounds_operators} evident.
\end{remark}

\medskip

Using these operators we define the following sequence of
functions:
\begin{equation}\label{def_nf_1}
f_n^{(0)} \isdef  1 \,, \quad  f_n^{(1)} \isdef  \mathcal M_n^{\rm
i}(1) \,, \quad f_n^{(2)} \isdef \mathcal M_n^{\rm e}
(f_n^{(1)})\,,
\end{equation}
and
\begin{equation}\label{def_nf_2}
f_n^{(2k+1)} \isdef  \mathcal M_n^{\rm i} (f_n^{(2k)}) \,, \quad
f_n^{(2k+2)} \isdef  \mathcal M_n^{\rm e} (f_n^{(2k+1)})\,, \quad
 k \in
\N\,.
\end{equation}
Observe that each of these formulas defines in fact two analytic
functions, one inside and one outside of the circle where the
corresponding integral is taken. In order to avoid cumbersome
notation we will not use the subindices $\rm i$ and $\rm e$ in this
case, but the reader should bear in mind this fact.

By \eqref{bounds_operators},
\begin{align*}
|z|\neq r \quad \Rightarrow \quad |f_n^{(1)}| & = |  \mathcal
M_n^{\rm i}(1)|\leq \frac{\Lambda r}{\tau
^2}\,  \frac{  r^n}{\left||z|-r \right|}\,, \\
|z|\neq 1/r \quad \Rightarrow \quad |f_n^{(2)}| & = |  \mathcal
M_n^{\rm e} (f_n^{(1)}) |\leq \frac{\Lambda \tau ^2}{r}\, \frac{
 r^n}{\left||z|-1/r \right|}\, \|f_n^{(1)}\|_{\T_{1/r}}\leq \frac{
\Lambda^2
r^{2n}}{\frac{1}{r}-r}\, \frac{ 1}{\left||z|-1/r \right|}\,,\\
|z|\neq r \quad \Rightarrow \quad |f_n^{(3)}| & = |  \mathcal
M_n^{\rm i} (f_n^{(2)}) |\leq  \frac{\Lambda r}{\tau ^2}\,  \frac{
 r^n}{\left||z|-r \right|}\, \|f_n^{(2)}\|_{\T_{r}}\leq \frac{r}{\tau ^2}\, \frac{ \Lambda^3
r^{3n}}{\left(\frac{1}{r}-r\right)^2}\, \frac{ 1}{\left||z|-r
\right|}\,,
\end{align*}
and in general, for $k\in \N$,
\begin{equation}\label{estimatesFk}
|f_n^{(k)}(z)| \leq \begin{cases} \dfrac{r}{\tau ^2}\, \dfrac{
\Lambda^k r^{kn}}{\left(\frac{1}{r}-r\right)^{k-1}}\, \dfrac{
1}{\left||z|-r
\right|}\,, & \text{if $k$ is odd, and }|z|\neq r\,, \\
 \dfrac{ \Lambda^k
r^{kn\strut}}{\left(\frac{1}{r}-r\right)^{k-1}}\, \dfrac{
1}{\left||z|-1/r \right|}\,, & \text{if $k$ is even, and }|z|\neq
1/r\,.
\end{cases}
\end{equation}

\begin{remark}
There are a number of equivalent definitions for the functions
$f_n^{(1)}$. For instance,
\begin{equation}\label{defFn}
-f_n^{(1)}(z)=  \frac{1}{2\pi i}\, \oint_{{\T_r}} \frac{t^n
\DD_{\rm i}^2(t)}{ w(t)(t-z)} \, dt=\frac{1}{2\pi i\, \tau^2}\,
\oint_{{\T_r}}
 \frac{t^n \FF(t) }{ t-z} \, dt=\frac{1}{2\pi i\, \tau^2}\,
\oint_{{\T_1}}
 \frac{t^n  \FF(t)}{ t-z} \, dt\,.
\end{equation}
We can rewrite the integrand directly in terms of the weight $w$
recalling that on ${\T_1}$,
$$
\frac{  \FF  (t)}{ \tau^2 }=\exp\left( \frac{1}{ \pi i}\, \dashint
_{{\T_1} } \frac{ \log w(s)\, ds}{s-t}\right)\,,
$$
where the integral is understood as its principal value. If
\begin{equation}\label{LaurentExpansion}
\FF(z)=\sum_{k=-\infty}^{+\infty } c_k z^k
\end{equation}
is the Laurent expansion of $\FF$ in the annulus $\rho <|z| <1/\rho
$ and $\mathcal P_+$ is the Riesz projection onto $H^2(\D_1)$, then
for $ z\in \D_r$,
\begin{equation}\label{projectionForF}
f_n^{(1)}(z)= -\frac{1}{\tau^2}\, \mathcal P_+ \left(z^n \FF(z)
\right)= -\frac{z^n}{\tau ^2}\,\times \left( \text{Laurent series of
$\FF$ truncated at $z^{-n}$}\right) \,,
\end{equation}
which establishes a connection of the $f_n^{(k)}$'s with Hankel
operators. It also shows that for $\rho <|z| < r$,
$$
\lim_{n \to \infty} \frac{f_n^{(1)}(z)}{z^n}=-\frac{\FF(z)}{\tau
^2}\,,
$$
and convergence is uniform on each compact subset of this annulus.

Finally, taking into account that
$$
f_n^{(1)}(0)=-\frac{1}{2\pi i\, \tau^2}\, \oint_{{\T_1}}
 t^{n-1}  \FF(t) \, dt=- \frac{c_{-n}}{\tau^2} \,,
$$
it is easy to see that sequence $\{ f_n^{(1)}\}$ satisfies the
following recurrence:
\begin{equation}\label{recForFn}
    f_{n+1}^{(1)}(z)=z f_{n}^{(1)} - \frac{c_{-n-1}}{\tau ^2} \,,  \quad z \in \D_1\,.
\end{equation}

%


\end{remark}

\medskip

Let
\begin{equation}\label{defSrepresentation}
\EE_n (z)    \isdef   \sum_{k=0}^\infty f_n^{(2k)}(z)\,, \quad
|z|\neq 1/r\,, \qquad \text{and} \qquad \II_n(z)  \isdef
\sum_{k=0}^\infty f_n^{(2k+1)}(z)\,, \quad |z|\neq r\,.
\end{equation}
Bounds \eqref{estimatesFk} show that the series in the right hand
sides of \eqref{defSrepresentation} converge absolutely and locally
uniformly in the respective regions of definition if
\begin{equation}\label{conditionConvergence}
 \Lambda  r^{ n} <\frac{1}{r}-r\quad \Longleftrightarrow \quad
 n>\frac{\log\left( \frac{1/r-r}{\Lambda}\right)}{\log r}\,.
\end{equation}
In this case there exists a constant $C>0$ depending on $r$ and
$\Lambda$ only, such that
\begin{equation}\label{uniformBoundsIandE}
|\EE_n (z)|\leq \frac{C}{||z|-1/r|}\,, \text{ for } |z|\neq 1/r\,,
\quad \text{and} \quad |\II_n (z)|\leq \frac{C}{||z|-r|}\,, \text{
for } |z|\neq r\,.
\end{equation}
In the sequel we use $C$ (eventually, with a subindex) to denote
some irrelevant constants, different in each appearance, whose
dependence or independence on the parameters will be stated
explicitly.

Now we may state one of the main results of the paper:
\begin{theorem}
\label{thm:nice case} Let $w$ be a strictly positive analytic
weight on the unit circle $\T_1$, the constant $\rho$ as defined
in \eqref{Def_R_bis}--\eqref{Def_R}, and constant $r$ with $\rho
<r<1$ fixed. Then with the notations introduced above, for every
$n\in \N$ satisfying \eqref{conditionConvergence} the following
formulas hold:
\begin{equation}\label{representation}
\Phi_n(z)=\begin{cases} \tau ^{-1} z^n D_{\rm e}(z) \,\EE_n
(z)\,, & \text{if } |z|> 1/r\,; \\
\tau ^{-1} z^n D_{\rm e}(z) \,\EE_n (z) -\dfrac{\tau\,\II_n(z)}{
D_{\rm
 i}(z)}\,, & \text{if } r < |z|< 1/r\,; \\
 -\dfrac{\tau\,\II_n(z)}{
D_{\rm
 i}(z)}\,, & \text{if }  |z|< r\,.
\end{cases}
\end{equation}
\end{theorem}
We emphasize that in each of the regions above $\EE_n$ and $\II_n$
have the meaning given in \eqref{defSrepresentation}, and are not
obtained in general by analytic continuation from one domain to
another.

\begin{remark}\label{remark:valid_everywhere}
The canonical formulae expressed in \eqref{representation} appear to
be valid for $z$ bounded away from the circles of radius $r$ and
$1/r$.  However, Theorem \ref{thm:nice case} holds true for any $r$
with $\rho < r < 1$, and so the formulae, in fact, are valid up to
the boundaries, by analytic continuation. Indeed, by considering
$\hat{r}$ with $r < \hat{r}<1$, the last formula in
\eqref{representation} is seen to remain valid by analytic
continuation, for $|z|<\hat{r}$. Similarly, the first formula is
valid by analytic continuation for $|z| > 1/\hat{r}$. For the middle
formula in \eqref{representation}, one considers $\tilde{r}$ with
$\rho < \tilde{r} < r$, which effectively extends the region where
the middle formula is valid. Therefore, the formulae contained in
\eqref{representation} are each in fact valid on the closure of the
domains stated.
\end{remark}

\begin{remark}
The second formula in \eqref{representation} expresses the
polynomial $\Phi_n$ on $\T_1$ as a sum of ``inner'' and ``outer''
terms. An alternative two-term expression for $\Phi_n$ on $\T_1$ can
be found in \cite[Formula (III.7)]{Geronimo/Case:79}. In the latter
formula we can identify the inner and outer Szeg\H{o} functions, but
the other factors are given implicitly in terms of the solutions of
equations from the inverse scattering theory.
\end{remark}

\medskip

Taking into account \eqref{tau}, the following
representation of the Verblunsky coefficients is straightforward: 
\begin{corollary}\label{cor:verblunsky}
With $w$ as in Theorem \ref{thm:nice case}, we have for each $n\in
\N$,
$$
\overline{\alpha _n}=\tau^2 \II_{n+1}(0)=\tau^2 \sum_{k=0}^\infty
f_{n+1}^{(2k+1)}(0)\,.
$$
\end{corollary}
Observe that under assumptions \eqref{conditionConvergence} the
series on the right hand side converges absolutely.

The Riemann-Hilbert asymptotic analysis that we carry out in
Section \ref{sec:RH_analysis} allows us to find also a
representation for the leading coefficients of the orthonormal
polynomials. For this purpose we define the ``shifted'' sequence
with respect to \eqref{def_nf_1}--\eqref{def_nf_2}:
\begin{align*}
g_n^{(0)} & \isdef 1\,, \quad g_n^{(1)} \isdef  \mathcal M_n^{\rm
e}(1) \,, \quad g_n^{(2)} \isdef \mathcal M_n^{\rm i}
(g_n^{(1)})\,, \\ \intertext{and} g_n^{(2k+1)} & \isdef \mathcal
M_n^{\rm e} (g_n^{(2k)}) \,, \quad g_n^{(2k+2)} \isdef \mathcal
M_n^{\rm i} (g_n^{(2k+1)})\,, \quad
 k \in
\N\,.
\end{align*}
Inequalities \eqref{bounds_operators} imply again that for $k \in
\N$,
\begin{equation}\label{estimatesGk}
|g_n^{(2k)}(0)| \leq \dfrac{ \Lambda^{2k}
r^{2kn-1}}{\left(\frac{1}{r}-r\right)^{2k-1}}\,.
\end{equation}
\begin{theorem}\label{prop:leading_coeff}
With $w$ as in Theorem \ref{thm:nice case}, and for $n$ satisfying
\eqref{conditionConvergence}, the leading coefficient $\kappa_{n}$
of $\varphi_n$ has the following representation:
$$
\kappa_{n}^2=\frac{\tau^2}{2 \pi}\,  \sum_{k=0}^\infty
g_n^{(2k)}(0) \,.
$$
\end{theorem}

\section{Asymptotic behavior} \label{sec:asymptotics}

By \eqref{estimatesFk} it is immediate that for $n$ satisfying
\eqref{conditionConvergence} and for $N=0, 1, 2,\dots $,
\begin{equation}\label{errorTruncation1}
    \left| \EE_n (z)  -   \sum_{k=0}^N f_n^{(2k)}(z)
     \right| \leq \frac{C_N}{\left| |z|-1/r\right|}\,  r^{(2N+2)
    n}\,, \quad |z|\neq 1/r\,,
\end{equation}
where the constant $C_N$ depends only on $r$, $N$ and $\Lambda$,
but neither on $n$ nor on $z$. Analogously,
\begin{equation}\label{errorTruncation2}
    \left| \II_n(z)  - \sum_{k=0}^N f_n^{(2k+1)}(z)
    \right| \leq \frac{C_N}{\left| |z|-r\right|}\,  r^{(2N+3)
    n}\,, \quad |z|\neq r\,,
\end{equation}
where $C_N$ has a similar meaning as above. Here a similar
observation concerning the validity of these bounds up to the
boundary, made in Remark \ref{remark:valid_everywhere}, applies.
These bounds show that \eqref{representation} allows us to obtain an
asymptotic expression for $\{\Phi_n\}$ of an arbitrarily high order.
Let us discuss the consequences of truncating $\II_n$ and $\EE_n$ in
\eqref{representation} at their first terms.

For instance, from Corollary \ref{cor:verblunsky} and Theorem
\ref{prop:leading_coeff} we get the following estimates of the
Verblunsky coefficients $\alpha _n$ and the leading coefficients
$\kappa _n$:
\begin{corollary}\label{cor:verblunsky2}
Let $w$ be a strictly positive analytic weight on the unit circle
$\T_1$. With the notation introduced above and for each $n\in \N$,
$$
\overline{\alpha _n}=- \frac{1}{2\pi i}\, \oint_{{\T_r}} t^{n} \,
\FF(t) \, dt +\mathcal O (r^{3n})=-c_{-n-1}+\mathcal O (r^{3n}) \,,
$$
where $c_{-n-1}$ is the corresponding Laurent coefficient of $\FF$
in \eqref{LaurentExpansion}. Consequently,
\begin{equation}\label{boundForVerbl}
|\alpha _n|\leq  \Lambda\, r^{n+1} +\mathcal O (r^{3n}) \,.
\end{equation}
Furthermore,
$$
\left| \kappa _n - \frac{\tau^2}{2 \pi} \right|\leq \dfrac{
\Lambda^{2} r^{2n}}{1-r^2}\,.
$$
\end{corollary}

Another straightforward consequence of truncating $\II_n$ and
$\EE_n$ at their first terms is the exterior asymptotics of $\{
\Phi_n\}$:
\begin{proposition} \label{cor:asymptotics_ouside}
Let $w$ be a strictly positive analytic weight on the unit circle
$\T_1$, and the constant $\rho$ as defined in
\eqref{Def_R_bis}--\eqref{Def_R}. Then
\begin{equation}\label{exterior_asymptotics}
\lim_{n \to \infty}\frac{\Phi_n(z)}{z^n}= \frac{ D_{\rm
e}(z)}{\tau}
\end{equation}
locally uniformly for $|z|> \rho $ and convergence takes place
with a geometric rate. In particular, for any compact set $K$ in
$\{z\in \C:\, |z|> \rho \text{ and } D_{\rm e}(z)\neq 0 \}$ there
exists $m=m(K)\in \N$ such that for $z\in K$, we have that
$\Phi_n(z)\neq 0$ for all $n \geq m$.

Furthermore, each zero $\zeta$ of $D_{\rm e}$ in $\rho <|z|<1$
attracts as many zeros of $\Phi_n$ as its multiplicity, which tend
to $\zeta$ with a geometric rate.
\end{proposition}

As mentioned in the Introduction, formula
\eqref{exterior_asymptotics} (or the equivalent formula
\eqref{ExteriorSzego}) shows that in the analytic case the
Szeg\H{o}'s exterior asymptotics can be continued analytically
through $\T_1$ up to the critical circle $\T_\rho $ (cf.\
\cite{Nevai/Totik:89} or \cite[\S 7.1]{Simon04a}). The rate of
convergence is a consequence of \eqref{representation}. For
instance, when $r<|z|<1/r$, by \eqref{estimatesFk} and
\eqref{errorTruncation1}--\eqref{errorTruncation2},
$$
\left| \frac{\Phi_n(z)}{z^n}-\frac{ D_{\rm
e}(z)}{\tau}\right|=\left|\frac{D_{\rm e}(z)}{\tau}\left(
\EE_n(z)-1\right) - \frac{\tau}{D_{\rm i}(z)}\, \II_n(z) \right|
\leq \frac{C}{(1/r-|z|) (|z|-r)}\, \left| \frac{r}{z}\right|^n\,,
$$
where $C$ depends on $r$ and $\Lambda$ only.

The second statement of the Proposition is a direct application of
this fact and of the Hurwitz's theorem.

Regarding the zeros of $\Phi_n$, we see that besides those
matching the zeros of $D_{\rm e}$ outside $\T_\rho $, the rest
accumulates on $\overline{\D}_\rho $. Nevai-Totik's theorem
\cite{Nevai/Totik:89} states that in fact $\T_\rho$ separates the
regions of uniformly bounded and unbounded number of zeros of
$\Phi_n$; more precisely,
$$
\rho =\inf \left\{ r>0:\, \exists \, N(r)<\infty \text{ so that }
\forall n, \Phi_n(z) \text{ has } \leq N(r) \text{ zeros in }
|z|>r \right\}
$$
(see also \cite[\S 7.1]{Simon04a}).  Furthermore, by
\eqref{representation},
\begin{equation}\label{roughEstimatePhi}
\Phi_n(z)=-\frac{\tau }{D_{\rm i}(z)}\left(f_n^{(1)}(z)+h_n^{(1)}(z)
\right)\,, \quad z\in \D_r\,,
\end{equation}
and from \eqref{errorTruncation1}--\eqref{errorTruncation2} it
follows that there exists a constant $C$, depending only on
$\Lambda$ and $r$, such that if $\rho <r'<r$, then
\begin{equation}\label{boundH}
\left|h_n^{(1)}(z)\right|\leq \frac{C}{r-r'}\, r^{3n}\,, \quad z \in
\D_{r' }\,.
\end{equation}
As pointed out above, in a generic case the vast majority of zeros
distributes uniformly at the critical circle $\T_\rho$, although
some of them can remain inside, and by
\eqref{roughEstimatePhi}--\eqref{boundH}, stay close to the zeros of
$f_n^{(1)}$ in $\D_\rho$.

Let us make these statements more precise. For that purpose we
must impose some additional conditions on our weight $w$ (or
function $\FF$).
\begin{proposition} \label{prop:residues}
Assume that there exists $0\leq \rho '<\rho $ such that $D_{\rm
e}$ can be continued to the exterior of the circle $\T_{\rho'}$,
as an analytic function whose only singularities are on the circle
$\T_\rho$, and these are all isolated. Denote by $a_1, \dots, a_u$
the singularities (whose number is finite) of $D_{\rm e}$ on
$\T_\rho $,
\begin{equation}\label{equal_poles}
|a_1|= \dots= |a_u|=\rho \,.
\end{equation}
Then for $\rho <r'<r$ there exist constants $0\leq \delta=\delta
(r')<1$ and $C=C(r')<+\infty$ such that for $\rho '<|z|\leq r'$
and $n\in \N$,
\begin{equation}\label{case_finite_number_sing2}
\left|\Phi_n(z)-  z^n\,  \frac{D_{\rm e}(z)}{\tau } -\frac{1}{\tau
\, D_{\rm i}(z)} \sum_{k=1 }^u \res_{t=a_k } \left( \FF (t) \,
\frac{t^n }{ t-z } \right)\right| \leq C  \left( \rho ^n \delta
^n+ r^{3n}\right)\,.
\end{equation}
Furthermore, for every compact set $K\subset\D_\rho$ there exist
constants $0\leq \delta=\delta (K)<1$ and $C=C(K)<+\infty$ such that
for $z \in K$,
\begin{equation}\label{asymptotics_finite_sing}
\left| \Phi_n(z)-  \frac{1}{\tau \, D_{\rm i}(z)} \sum_{k=1 }^u
\res_{t=a_k } \left( \FF (t)   \, \frac{t^n }{ t-z }
\right)\right| \leq C  \left( \rho ^n \delta ^n+ \rho
^{3n}(1+\varepsilon )^{3n}\right)\,.
\end{equation}
\end{proposition}
\begin{remark}
If $\FF$ can be continued as an analytic function with a finite
number of isolated singularities to whole disc $\D_1$, then we may
take $\delta =0$ in the right hand sides in
\eqref{case_finite_number_sing2}--\eqref{asymptotics_finite_sing}.
Otherwise the right hand side in \eqref{asymptotics_finite_sing} may
be replaced by an estimate of the form $C \rho ^n \delta ^n$.
\end{remark}


\begin{remark}
This Proposition tells us that in general all the relevant
information for the asymptotics of $\Phi_n$'s in $\D_r$ comes from
the singularities of the exterior Szeg\H{o} function $D_{\rm e}$
on $\T_\rho $ (that is, from the first singularities of $D_{\rm
e}$ we meet continuing it analytically inside the unit disc).
\end{remark}


\begin{remark}
In order to isolate the zeros of $\Phi_{n}$, one must be able to
analyze the approximation to $\Phi_{n}$ afforded by
\eqref{case_finite_number_sing2} and
\eqref{asymptotics_finite_sing}. For example, zero-free regions
may be determined by \emph{(i)} establishing zero-free regions for
the approximation, and \emph{(ii)} bounding $\Phi_{n}$ away from
zero using the error estimates. Similarly, isolating the zeros can
be done by first isolating the zeros of the approximation, and
then using a Rouche' type argument for $\Phi_{n}$.
\end{remark}

\begin{remark} \label{remark:vicinity} It is interesting to note that Theorem \ref{thm:nice case} provides
global asymptotics for $\Phi_{n}$, and as a consequence,
Proposition \ref{prop:residues} above provides an asymptotic
description which is valid even in vicinities of the poles
$a_{k}$.
\end{remark}

\begin{remark} It is an interesting challenge to find general
assumptions concerning the measure of orthogonality that lead to an
effective characterization of the asymptotic behavior of the zeros
of $\Phi_{n}$.  Simon \cite{Simon04b, Simon04c} has placed
assumptions on the Verblunsky coefficients which imply that all
$a_k$'s are simple poles. In this case the approximation afforded by
\eqref{case_finite_number_sing2} and \eqref{asymptotics_finite_sing}
is completely explicit, and it is then straightforward (albeit
technically complicated) to isolate and provide a complete
asymptotic description for each and every zero of $\Phi_{n}$.  While
we could carry out this extension of \cite{Simon04b, Simon04c}, we
rather choose to investigate a more general class of measures of
orthogonality than considered in \cite{Simon04c}, and demonstrate
how it is possible to obtain the same detailed level of asymptotic
information concerning the zeros of $\Phi_{n}$. The first result in
this direction is the following Theorem.
\end{remark}

\begin{definition} 
Let $a$, $|a|<1$, be a pole of a function $f$ analytic in $|z|<1$.
We denote by $ \mult_{z=a}f(z)$ its multiplicity and say that $a$
is a \emph{dominant} pole of $f$ if for any other singularity $b$
of $f$, either $|a|>|b|$ or $|a|=|b|$, but then $b$ is also a pole
and
$$
\mult_{z=a}f(z) \geq \mult_{z=b}f(z)\,.
$$
\end{definition}
In the sequel we use the following notation: for $a\in \C$ and
$\varepsilon >0$,
$$
B_\varepsilon (a)\isdef \{z\in \C:\, |z-a|<\varepsilon  \}\,.
$$

\begin{theorem} \label{prop:many dominant}
Assume that there exists \/ $0\leq \rho '<\rho $ such that $D_{\rm
e}$ can be continued to the exterior of the circle $\T_{\rho'}$,
as a meromorphic function whose only singularities are on the
circle $T_\rho$. Denote by $a_1, \dots, a_u$ the poles (whose
number is finite) of $D_{\rm e}$ on $\T_\rho $, and assume that
the dominant poles of $D_{\rm e}$ are $a_1, \dots, a_\ell$, $\ell
\leq u$, and their multiplicity is $m$.

Let $\varepsilon >0$. Then for $\rho '<|z|\leq r-\varepsilon$, $z
\notin \cup_{k=1}^u B_\varepsilon (a_k)$, and $n\in \N$,
\begin{equation}\label{case_finite_number_dominant_sing}
\Phi_n(z)=  \frac{D_{\rm e}(z) }{\tau }\, z^n +\frac{1}{\tau
D_{\rm i}(z)} \, \sum_{k=1 }^\ell \binom{n}{m-1}\, a_k^{n-m+1}\,
\frac{D_{\rm i} (a_k)\, \widehat D_{\rm e} (a_k ) }{ a_k-z
 } \,   + h_n(z)   \,,
\end{equation}
where $ \widehat D_{\rm e} (a_k)=\lim_{z\to a_k} D_{\rm e}
(z)(z-a_k)^m $, $ k=1, \dots, \ell$. There exist a constant
$0<C<+\infty$ independent of $n$ and $\varepsilon $, and a
constant $0<\delta=\delta (\varepsilon )<1$, such that
\begin{equation}\label{newBOundForH}
    |h_n(z)|\leq \begin{cases}
    C \left( \rho ^n \delta ^n+ r^{3n}\right), & \text{if } \, m=1\,, \\
    \dfrac{C}{ \varepsilon ^{m-1}}\, n^{m-2} \, \rho ^n\,, & \text{if }\, m\geq 2\,.
    \end{cases}
\end{equation}

Furthermore, for every compact set $K\subset\D_\rho$ there exists a
constant $C=C(K) <\infty$ such that for $z \in K$, and $n\in \N$,
\begin{equation}\label{asymptotics_dominant_poles}
\left| \frac{\tau \, D_{\rm i}(z)}{a_1 ^{n-m+1}}\, \binom{n }{ m-1
}^{-1}\,\Phi_n(z)-
   \sum_{k=1 }^\ell \frac{D_{\rm i} (a_k)\,  \widehat D_{\rm e} (a_k) }{  a_k-z
 }   \, e^{2\pi i (n-m+1) \theta_k}\right|\leq \begin{cases}
C   \delta ^n, & \text{if } \, m=1\,, \\
    \dfrac{C}{   n} \,, & \text{if }\, m\geq 2\,,
 \end{cases}
\end{equation}
where
\begin{equation}\label{defTheta}
\theta_1=1, \quad \text{and} \quad
\theta_k=\frac{1}{2\pi}\,\left(\arg a_k - \arg a_1 \right),
  \qquad  k=2, \dots, \ell\, .
\end{equation}
In particular, on every compact set $K\subset \D_\rho $, for all
sufficiently large $n$ polynomials $\Phi_n$ can have at most $\ell
-1$ zeros, counting their multiplicities.
\end{theorem}

\begin{remark}
In the spirit of Proposition \ref{prop:residues}, the asymptotics
of $\Phi_n$'s is governed by the dominant poles of $\FF$ only.
However, the error term $h_n$ in
\eqref{asymptotics_dominant_poles} is small only away from all
poles $a_k$.

Observe that this result is applicable to weights of the form
$w(z)=|R(z) S(z)|^2$, $z\in \T_1$, where $R$ is a rational function
with at least one zero on $\T_\rho $ (or one pole on $\T_{1/\rho
}$), and $S$ is any function holomorphic and $\neq 0$ in any
annulus, containing $\{ \rho \leq |z| \leq 1/\rho \}$.
\end{remark}

\begin{remark}  As mentioned in Remark \ref{remark:vicinity} above, we may also obtain an
asymptotic description of $\Phi_{n}$ in a vicinity of each of the
poles $a_{k}$, and describe the behavior of the zeros near these
singular points.  We leave this exercise to the interested reader.

On the other hand, to investigate new phenomena, we consider in
Section \ref{sec:examples} the situation that there is a single
isolated singularity of $D_{\rm e}$ on the circle of radius $\rho$
which is an essential singularity.  For this case, we do carry out
the detailed analysis of $\Phi_{n}$ near the location of the
essential singularity.
\end{remark}

\begin{corollary} \label{cor:verblunskyFiniteCase}
Under assumptions of Theorem \ref{prop:many dominant}, the
Verblunsky coefficients satisfy
\begin{equation}\label{Verblunsky_asymptotics_dominant_poles}
 \alpha _n=
   -\sum_{k=1 }^\ell \binom{n+1}{m-1}\, \overline{a_k^{n-m+1 }\,
 D_{\rm i} (a_k)\, \widehat D_{\rm e} (a_k ) }  + \begin{cases}
    \mathcal O \left( \rho ^n \delta ^n\right), & \text{if } \, m=1\,, \\
    \mathcal O\left(  n^{m-2} \, \rho ^n\right)\,, & \text{if }\, m\geq 2\,.
    \end{cases}
\end{equation}
\end{corollary}
\begin{remark}
The case of Verblunsky coefficients with the asymptotics
\eqref{Verblunsky_asymptotics_dominant_poles} and $m=1$ was analyzed
in \cite{Simon04c}. The subcase $m=1$, $u=1$ has been thoroughly
studied in \cite{Simon04b}.
\end{remark}

\medskip

Let us discuss now how Theorem \ref{prop:many dominant} reveals the
behavior of the zeros of the polynomials in $\D_r$. We distinguish
two cases: the zeros approaching the critical circle $\T_\rho $ and
those staying inside $\D_\rho $. Under assumptions of Theorem
\ref{prop:many dominant} let us denote
\begin{equation}\label{G_n}
    G_n(z)\isdef \sum_{k=1 }^\ell
\frac{D_{\rm i} (a_k)\, \widehat D_{\rm e} (a_k ) }{ a_k-z
 }\, e^{2\pi i (n-m+1) \theta_k}\,,
\end{equation}
$\theta_k$ defined in \eqref{defTheta}. Numerical experiments show
that that the zeros of $\Phi_n$ approaching $\T_\rho $ distribute
radially uniformly, with gaps (``missing zeros'') at the poles
$a_k$ and at possible zeros of $G_n$ on $\T_\rho $. We will use
the following notation: for a given analytic function $f$ and for
$\varepsilon >0$,
\begin{equation}\label{nullSet}
    \mathcal Z(f) \isdef \{z:\, f(z)=0 \}\,, \quad \mathcal Z_\varepsilon (f )
    \isdef   \bigcup_{z\in \mathcal Z(f)} B_\varepsilon
    (z)\,.
\end{equation}
Also let
$$
\mathcal B(\varepsilon )\isdef \left\{z:\, \rho' + \varepsilon
<|z|<r- \varepsilon, \; z\notin \cup_{k=1}^u B_\varepsilon
(a_k)\cup \mathcal Z_\varepsilon (D_{\rm e} ) \right\} \,.
$$

\begin{theorem}
\label{thm:clock} Let $\varepsilon
>0$, and assume that the hypotheses of Theorem
\ref{prop:many dominant} are satisfied. Put
\begin{equation}\label{defVn}
 v_n \isdef \rho \,\binom{n}{m-1}^{1/n}\,.
\end{equation}
Then there is a $\delta
>0$ such that for all sufficiently large $n\in\N$ every pie slice
of the form
\begin{equation*}
\begin{split}
 \left\{z:\, v_n(1-\delta/2)
<|z|<v_n(1 + \delta/2) , \; \frac{\alpha +2 k_1
\pi}{n}<\arg\left(z \right)<\frac{\alpha  +2 k_2 \pi}{n}
\right\}\subset \mathcal B(\varepsilon )\setminus \mathcal
Z_\varepsilon(G_n )\,,
\end{split}
\end{equation*}
with appropriately chosen $\alpha \in \R$, contains exactly
$k=k_2-k_1 $ zeros of $\Phi_n$, $z^{(n)}_1,\dots, z^{(n)}_k$,
satisfying
\begin{equation}\label{asymptModEquidistr}
|z^{(n)}_i|=\rho \left(1+\frac{1}{n}\, \log
\binom{n}{m-1}+\mathcal O\left(\frac{1}{n}\right) \right)\,,
\end{equation}
and
\begin{equation}\label{asymptArgEquidistr}
\arg\left(z^{(n)}_{i+j}\right)-\arg\left(z^{(n)}_i\right)=\frac{2\pi
j}{n}+\mathcal O \left( \frac{1}{n^2} \right)\,.
\end{equation}

\end{theorem}

\begin{figure}[htb]
\centering \begin{overpic}[scale=0.95]{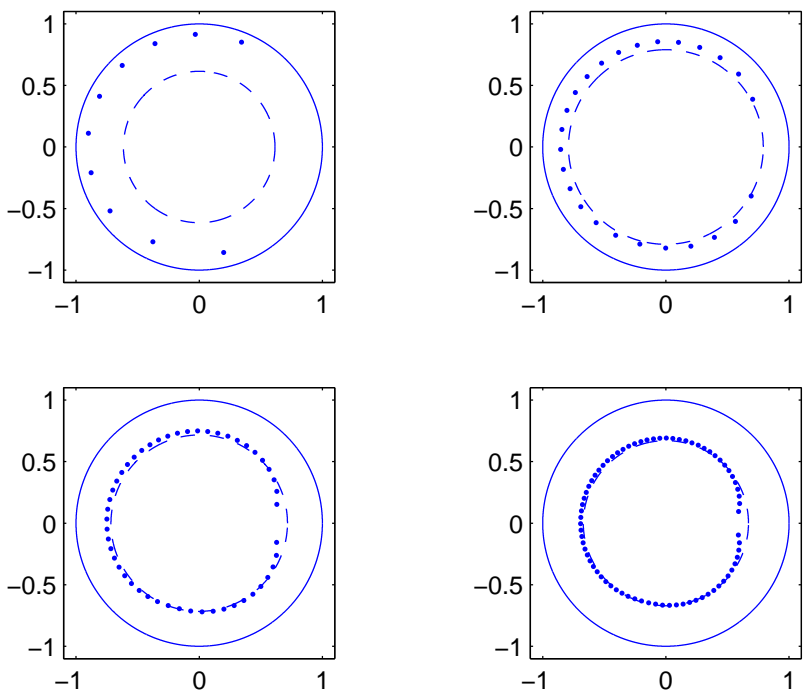}
\end{overpic}
\caption{Zeros of $\Phi_n$ for $n=10, 25, 50, 75$  (from left to
right and from top to bottom) with
$w(z)=|(z-1/2)^{10}(z-\zeta_1/2)^3|^2$, $z\in \T_1$, where
$\zeta_1=\exp(\pi i \sqrt{2})$. For comparison with the prediction
\eqref{asymptModEquidistr} we plot in each case with a dashed line
the circle of radius $\frac{1}{2} \left( 1+\frac{1}{n}\, \log
\binom{n}{m-1}\right)$.}\label{fig:zerosMultiplePoles2bisCC}
\end{figure}

\begin{remark}
We may consider an example for the situation described in Theorem
\ref{prop:many dominant}, taking
$w(z)=|(z-1/2)^{10}(z-\zeta_1/2)^3|^2$, $z\in \T_1$, with
$\zeta_1=\exp(\pi i \sqrt{2})$. The results of the numerical
experiments appear in Fig. \ref{fig:zerosMultiplePoles2bisCC}.
Observe how the zeros exhibit the equidistribution pattern close
to the circle of radius $0.5 \left( 1+\frac{1}{n}\, \log
\binom{n}{m-1}\right)$ predicted by Theorem \ref{thm:clock}, and
basically ``ignoring'' the non-dominant pole $\zeta _1/2$ of
$\FF$. An imaginative reader might observe though that the the
influence of the non-dominant pole is not totally negligible: Fig.
\ref{fig:zerosMultiplePoles2bisCC} shows a slight displacement of
the zeros of $\Phi_n$ towards the origin in a neighborhood of
$\zeta _1/2$.
\end{remark}

We can also be more precise about the accumulation set of zeros of
$\Phi_n$'s on compact sets of $\D_\rho $. Let
\begin{equation}\label{def:accumulationSet}
    Z \isdef \bigcap_{k\geq 1} \overline{\bigcup_{n\geq k} \mathcal Z( \Phi_n) }\,.
\end{equation}
Asymptotics \eqref{asymptotics_dominant_poles} shows that the
structure of $Z\cap \D_\rho $ depends on the relative positions of
$a_1, \dots, a_\ell$ on $\T_\rho $. Without loss of generality we
may assume  that $\theta_1=1, \theta_2, \dots, \theta_v$ is the
maximal subset of $\{ \theta_1, \dots, \theta_\ell\}$ linearly
independent over the rational numbers $\Q$. Then there exist unique
$r_{kj}\in \Q$, $k=1, \dots, \ell$, $j=1, \dots, v$, such that
$$
\theta_k=\sum_{j=1}^v r_{kj}\, \theta_j\,, \quad \ k=1, \dots,
\ell\,,
$$
(obviously, $r_{kj}=\delta_{kj}$, for $1\leq k \leq v$).
Furthermore, we may always take $r_{k1}\geq 0$, $k=1, \dots, \ell$.
Denote
$$
r_{k1}=\frac{p_k}{q_k}, \quad 0\leq p_k\leq q_k\in \Z\,, \quad k=1,
\dots, \ell\,.
$$
\begin{proposition}
\label{prop:szabados} Let $t\in Z\cap \D_\rho $. If $v=1$ (that is,
if all $\theta_k\in \Q$), then there exist $0\leq s_k<q_k$, $s_k\in
\Z$, $k=1,\dots,\ell$ such that
\begin{equation}\label{description_Szabados}
\sum_{k=1}^\ell \frac{D_{\rm i} (a_k)\,  \widehat D_{\rm e} (a_k) }{
a_k-t  }\, e^{2\pi i  s_k/q_{k} }=0\,.
\end{equation}
In this case $Z\cap \D_\rho$ is a finite set.

If $v\geq 2$, then additionally there exist $X_2, \dots, X_v\in \R$
such that
\begin{equation}\label{description_Szabados2}
\sum_{k=1}^\ell \frac{D_{\rm i} (a_k)\,  \widehat D_{\rm e} (a_k) }{
a_k-t  }\, e^{2\pi i \left( s_k/q_{k}+\sum_{j=2}^v r_{kj} \, X
_j\right)}=0\,.
\end{equation}
\end{proposition}

\begin{figure}
\centering
\begin{tabular}{ll}
\hspace{-1.5cm}\mbox{\includegraphics[scale=0.65]{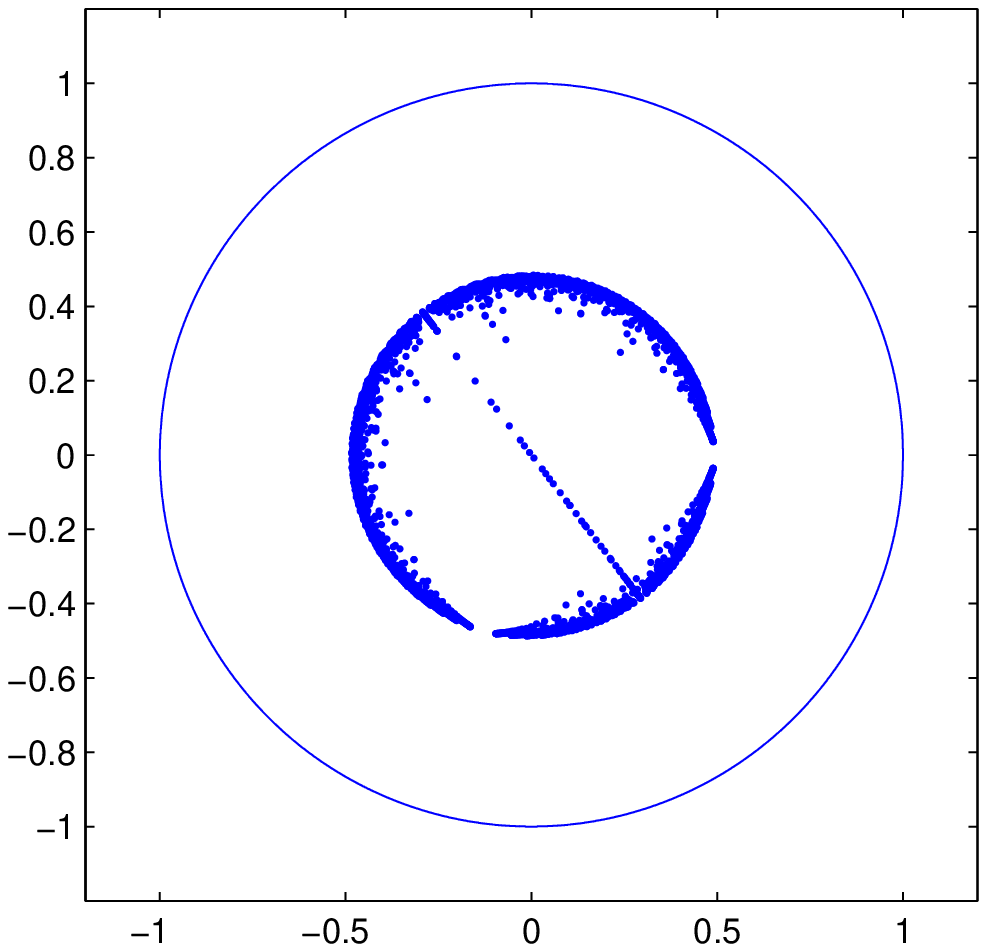}} &
\hspace{-2cm}\mbox{\includegraphics[scale=0.65]{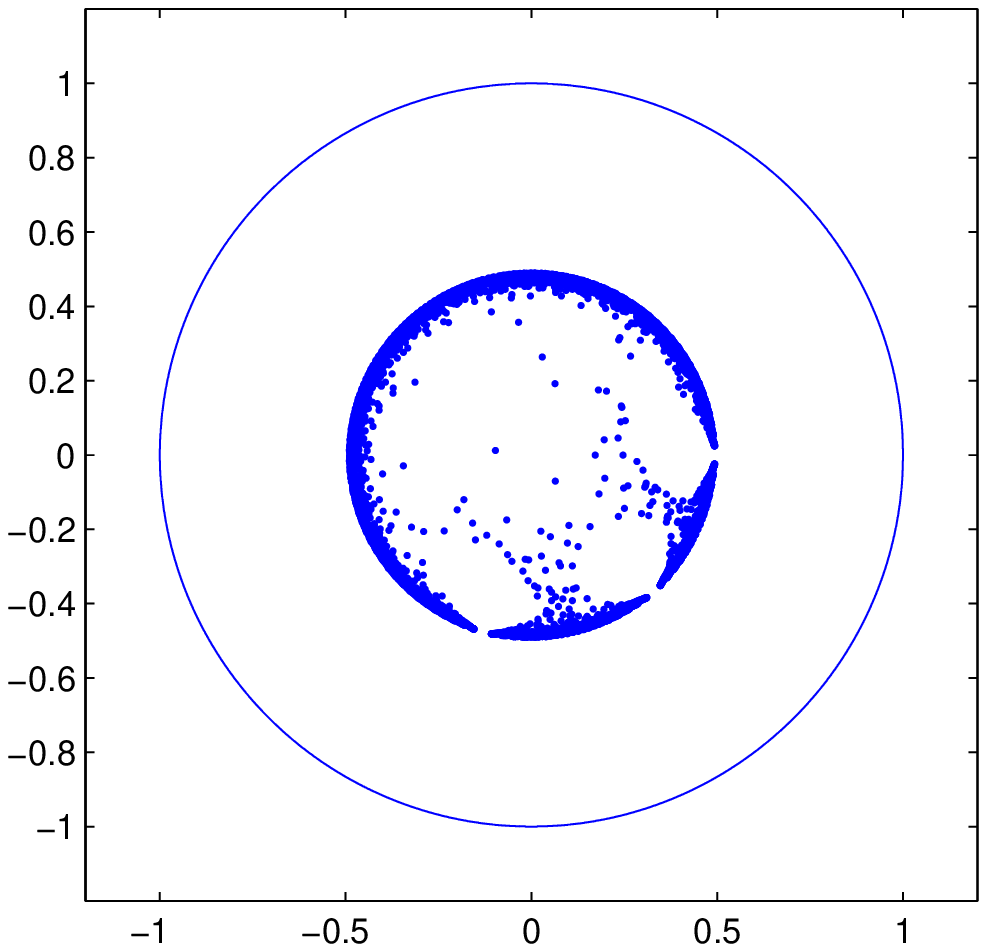}}
\end{tabular}
\caption{Zeros of $\Phi_n$ for $n=5, \dots, 75$ with
$w(z)=|(z-1/2)(z-\zeta_1/2)|^2$ (left) and
$w(z)=|(z-1/2)(z-\zeta_1/2)(z-\zeta_2/2)|^2$ (right), $z\in \T_1$,
where $\zeta_1=\exp(\pi i \sqrt{2})$, $\zeta_2=\exp(\pi i
\sqrt{3})$.}\label{fig:twoandthreepoles}
\end{figure}

\begin{remark} \label{remark:cirves}
As it was observed by Szabados \cite{Szabados79} (in the polynomial
case), for $v=2$ equation \eqref{description_Szabados2} describes an
algebraic curve of degree $\leq \ell$. In particular, if $\ell=2$,
$Z\cap \D_\rho$ is a circular arc or a diameter in $\D_\rho $ (see
Figure \ref{fig:twoandthreepoles}, left). For $v> 2$ equation
\eqref{description_Szabados2} describes a two-dimensional domain
bounded by algebraic curves (cf.\ Figure \ref{fig:twoandthreepoles},
right). Observe that we plot \emph{all} the zeros of $\Phi_n$ for
$n=5,\dots, 75$, in order to reveal the structure of $Z$ inside
$\T_{1/2}$.
\end{remark}

We can be more specific in the case of a \emph{rational} weight,
i.e.\ when $w$ has the following form:
\begin{equation}\label{rational_case}
    w(z)=\left|r(z) \right|^2\,, \quad |z|=1\,,
\end{equation}
where $r$ is a rational function non vanishing on $\T_1$; we may
write
$$
r(z)=  \frac{p_{\rm i} }{q_{\rm i}}(z)\, \frac{p_{\rm e} }{q_{\rm
e}}(z)\,,
$$
where $p_{\rm i}$ and $q_{\rm i}$ are polynomials without zeros in
$|z|\leq 1$, $p_{\rm e}$ and $q_{\rm e}$ are polynomials without
zeros in $|z|\geq 1$, and $p_{\rm e}(0)\neq 0$ and $q_{\rm
e}(0)\neq 0$.
Then,
$$
|r(z)|=\left|\frac{p_{\rm i} }{q_{\rm i}}(z)\, \frac{p_{\rm e}
}{q_{\rm e}}(z) \right|=\left|e^{i\alpha } \frac{p_{\rm i}
}{q_{\rm i}}(z)\, \frac{p_{\rm e} }{q_{\rm e}}(z) \, \frac{p_{\rm
i}^* }{p_{\rm i} }(z)\, \frac{q_{\rm i} }{q_{\rm i}^*}(z)
\right|=\left|e^{i\alpha } \frac{p_{\rm i}^* }{q_{\rm i}^*}(z)\,
\frac{p_{\rm e} }{q_{\rm e}}(z)  \right|=\left| R(z)\right|\,,
\quad \alpha \in \R\,,
$$
where $R$ has neither zeros nor poles in $|z|\geq 1$; since in
this identity we may multiply $R$ by any power of $z$ and any
constant with absolute value equal to 1, without loss of
generality we may assume that $0<R(\infty) < \infty$.
\begin{definition} \label{def:normalrepresentation}
We say that $w(z)=|R(z)|^2$, $|z|=1$ is a \emph{normalized
representation} of the rational weight $w$ if $R$ is a rational
function with neither zeros nor poles in $|z|\geq 1$ and such that
$0<R(\infty) < \infty$.
\end{definition}
If $w(z)=|R(z)|^2$, $|z|=1$, is a normalized representation of
$w$, by rewriting this identity on $ \T_1$ as
$$
w(z)=\left|R(z)\right|^2=R(z) \overline{R(z)}= R(z)
\overline{R}(1/z)\,,
$$
we immediately obtain that in this case,
\begin{equation}\label{szego_in_rational}
D_{\rm i}(z)=\overline{R}(1/z)\,, \quad D_{\rm e}(z)=
\frac{1}{R(z)}\,, \quad \tau= \frac{1}{R(\infty)}\,.
\end{equation}
Furthermore,
\begin{equation}\label{blaschke product}
\FF(z)=\frac{\overline{R}(1/z)}{R(z)}
\end{equation}
is a Blaschke product with its poles in $|z|< 1$ coinciding with
the zeros of $R$, and with its zeros in $|z|< 1$ matching the
poles of the same functions.

Let $a_1, \dots, a_\ell$ be the zeros of $R$ in $|z|< 1$. Now by
\eqref{representation},
\begin{equation}\label{caseRationalFunction}
\Phi_n(z)=   z^n\, \dfrac{R(\infty)}{R(z)}
+\frac{R(\infty)}{\overline{R}(1/z)}\, \sum_{k=1}^\ell \res_{t=a_k
} \left( \dfrac{\overline{R}(1/t)}{R(t)} \, \frac{t^n  }{  t-z }
\right)+h_n(z) \,, \quad |z|<1\,,
\end{equation}
where by \eqref{errorTruncation1}--\eqref{errorTruncation2},
$|h_n(z)|\leq C/(r-|z|)\, r^{3n}$, being $C$ a constant depending
on $\Lambda$ and $r$ only.

In particular, if all the zeros of $R$ in $|z|<1$ are simple, we
get
\begin{equation}\label{caseRationalFunction_onepole}
\Phi_n(z)=   z^n\, \dfrac{R(\infty)}{R(z)}
+\frac{R(\infty)}{\overline{R}(1/z)} \, \sum_{k=1}^\ell
\dfrac{\overline{R}(1/a_k)}{R'(a_k)} \, \frac{a_k^n  }{ a_k-z }
+h_n(z) \,, \quad |z|<1\,.
\end{equation}
\begin{remark}
As mentioned in the introduction, Szabados \cite{Szabados79}
obtained a formula similar to \eqref{caseRationalFunction_onepole}
for the case of a polynomial weight; he also observed the
equidistribution of the zeros approaching the critical circle
$\T_\rho $. Rational weights with one dominant singularity have been
discussed also by Ismail and Ruedemann in \cite{Ismail92}.
\end{remark}
\begin{remark}
To contrast some of these results with a case of a more severe
singularity, in Section \ref{sec:examples} we consider the weight of
orthogonality
\begin{equation}\label{weightExampleSing}
w(t)=w_1(t)=\left|\exp\left( \frac{1}{a-t}\right) \right|^2\,, \quad
t \in \T_1\,,
\end{equation}
with $0<a<1$, as well as its inverse, $1/w$. For $w$ in
\eqref{weightExampleSing} the asymptotics of the Verblunsky
coefficients is
$$
\alpha_n=-\frac{1}{2\sqrt{\pi}}\, t_+^{n} \FF(t_+)\, \left(
\frac{a}{n}\right)^{3/4} \left( 1+\mathcal O \left(
\frac{1}{n^{1/2}}\right)\right)\,, \quad n \to \infty\,,
$$
(see Eq.\ \eqref{verblunskyForEssential}), where $t_+$ is the
solution of the equation
$$
\frac{1}{t}-\frac{1}{n+1} \, \left(\frac{t}{a t-1}+\frac{1}{t-a}
\right)=0
$$
satisfying
\begin{equation*}
    t_+= a + \sqrt{\frac{a}{n+1}}+\mathcal O\left(\frac{1}{n}\right)\,, \quad n \to
    \infty\,,
\end{equation*}
where we take the positive square root.

For the weight $1/w$, with $w$ defined in \eqref{weightExampleSing},
the asymptotics of the Verblunsky coefficients is
$$
\alpha_n=-\frac{1}{\sqrt{\pi}}\, \left( \frac{a}{n}\right)^{3/4}
 \, \Re \left( t_+^{n} \FF(t_+)\right) \, \left( 1+\mathcal O \left(
\frac{1}{n^{1/2}}\right)\right)\,, \quad n \to \infty\,,
$$
(see Eq.\ \eqref{verblunskyForEssential}), but $t_+$ is now the
solution of the equation
$$
\frac{1}{t}+\frac{1}{n+1} \, \left(\frac{t}{a t-1}+\frac{1}{t-a}
\right)=0
$$
satisfying
\begin{equation*}
    t_+= a + i\,\sqrt{\frac{a}{n+1}}+\mathcal O\left(\frac{1}{n}\right)\,, \quad n \to
    \infty\,,
\end{equation*}
where we take again the positive square root.

\begin{figure}[htb]
\centering
\begin{tabular}{ll}
\hspace{-1.5cm}\mbox{\includegraphics[scale=0.65]{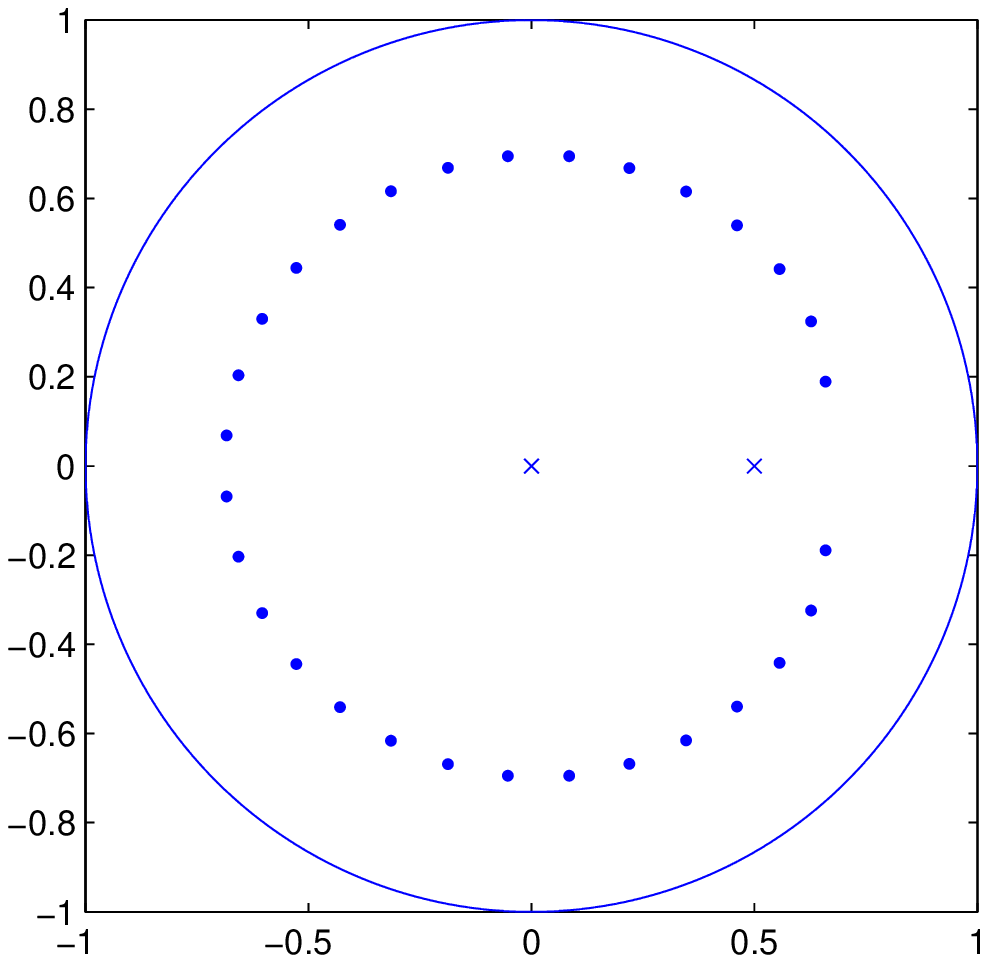}}
&
\hspace{-2cm}\mbox{\includegraphics[scale=0.65]{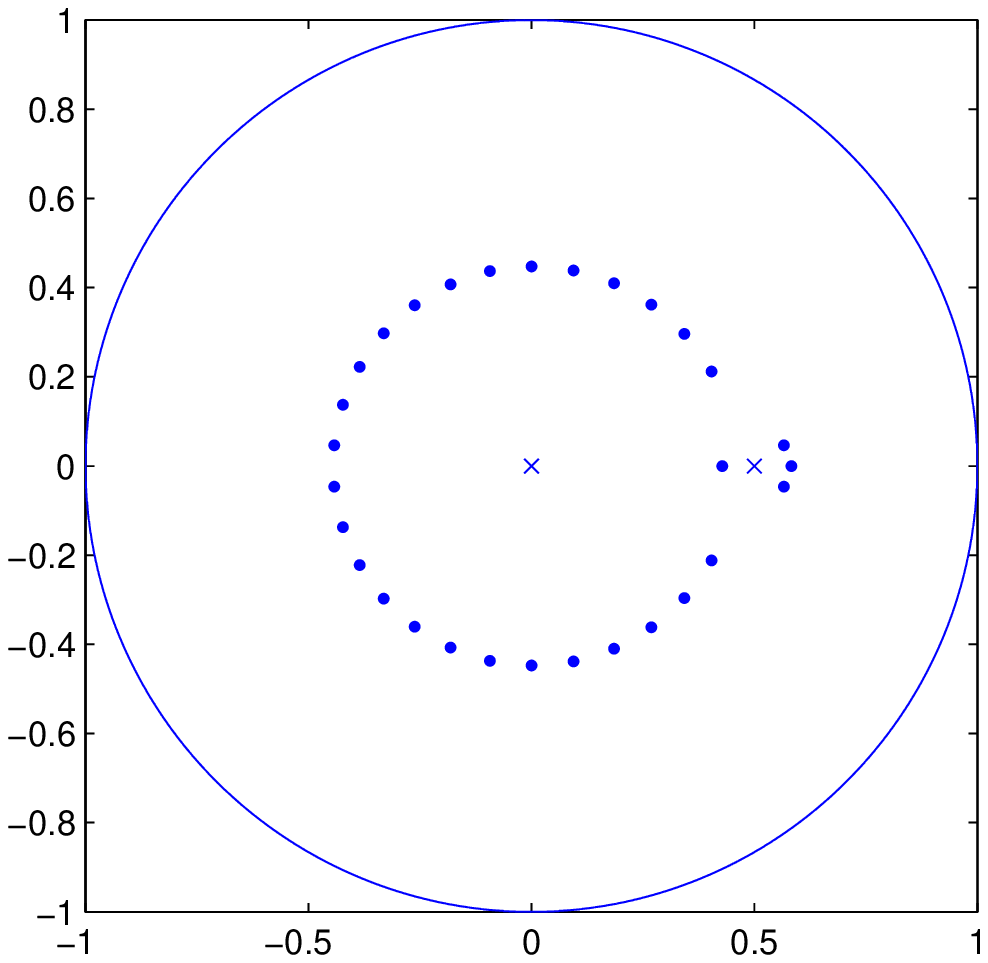}}
\end{tabular}
 \caption{Zeros of $\Phi_{30}$ for weights $w$ (left) and $1/w$ (right), with $w$
given in \eqref{weightExampleSing} and $a=1/2$.
}\label{fig:Essential30_dominantExample}
\end{figure}

Also the zero behavior in both cases is different. In order to
entice the reader we show in Fig.\
\ref{fig:Essential30_dominantExample} images of the zeros of
$\Phi_n$ (for $n=30$ and $a=1/2$), whose complete asymptotic
description appears in Section \ref{sec:examples}.
\end{remark}

\section{Steepest descent analysis for orthogonal polynomials}
\label{sec:RH_analysis}

Observe that from the orthogonality
conditions \eqref{orthogonalityConditions} it follows that
$$
\oint_{\T_1} \varphi_n(z)  \overline{z^k}\, w(z)
|dz|=\frac{1}{i}\, \oint_{\T_1} \varphi_n(z)  z^{n-k-1} \,
\frac{w(z)}{z^n} dz=0\quad \text{for } k=0, 1, \dots, n-1\,.
$$
The reversed polynomials $\varphi^*(z)=z^n
\overline{\varphi_n}(1/z)$ also satisfy the following
orthogonality conditions:
$$
 \oint_{\T_1} \varphi^*_{n-1}(z)  z^{k} \,
\frac{w(z)}{z^n} \, dz=i\, \oint_{\T_1}
\overline{\varphi_{n-1}(z)} z^k\, w(z) |dz|=\begin{cases} 0, &
k=0, 1, \dots, n-2, \\ i/\kappa _{n-1}, & k=n-1\,.
\end{cases}
$$
Hence,
$$
Y(z)=\begin{pmatrix} \Phi_n(z) & \displaystyle \dfrac{1}{2\pi i}\,\oint_{\strut\T_1}
\dfrac{\Phi_n(t) w(t)\, dt}{t^n(t-z)}  \\
-2\pi \kappa _{n-1} \varphi^*_{n-1}(z) & - \, \displaystyle
\dfrac{\kappa _{n-1}}{ i}\,\oint_{\strut\T_1}
\dfrac{\varphi^*_{n-1}(t) w(t)\, dt}{t^n(t-z)}
\end{pmatrix}
$$
is a unique solution of the following Riemann-Hilbert problem: $Y$
is holomorphic in $\C\setminus {\T_1}$,
$$
Y_+(t)=Y_-(t)\, \begin{pmatrix} 1 & w(t)/t^n \\ 0 & 1
\end{pmatrix}\,, \quad z\in \T_1\,,
$$
and
$$
\lim_{z\to \infty} Y(z)\,\begin{pmatrix} z^{-n} & 0 \\ 0 & z^n
\end{pmatrix}=I\,,
$$
where $I$ is the $2 \times 2$ identity matrix. Observe also that
\begin{equation}\label{kappa}
\kappa_{n-1}^2=-\frac{Y_{21}(0)}{2\pi}\,,
\end{equation}
(see \cite{MR2000e:05006} or \cite{MR2000g:47048} for more details
on the Riemann-Hilbert characterization of orthogonal
polynomials).

\begin{figure}[htb]
\centering \begin{overpic}[scale=0.65]{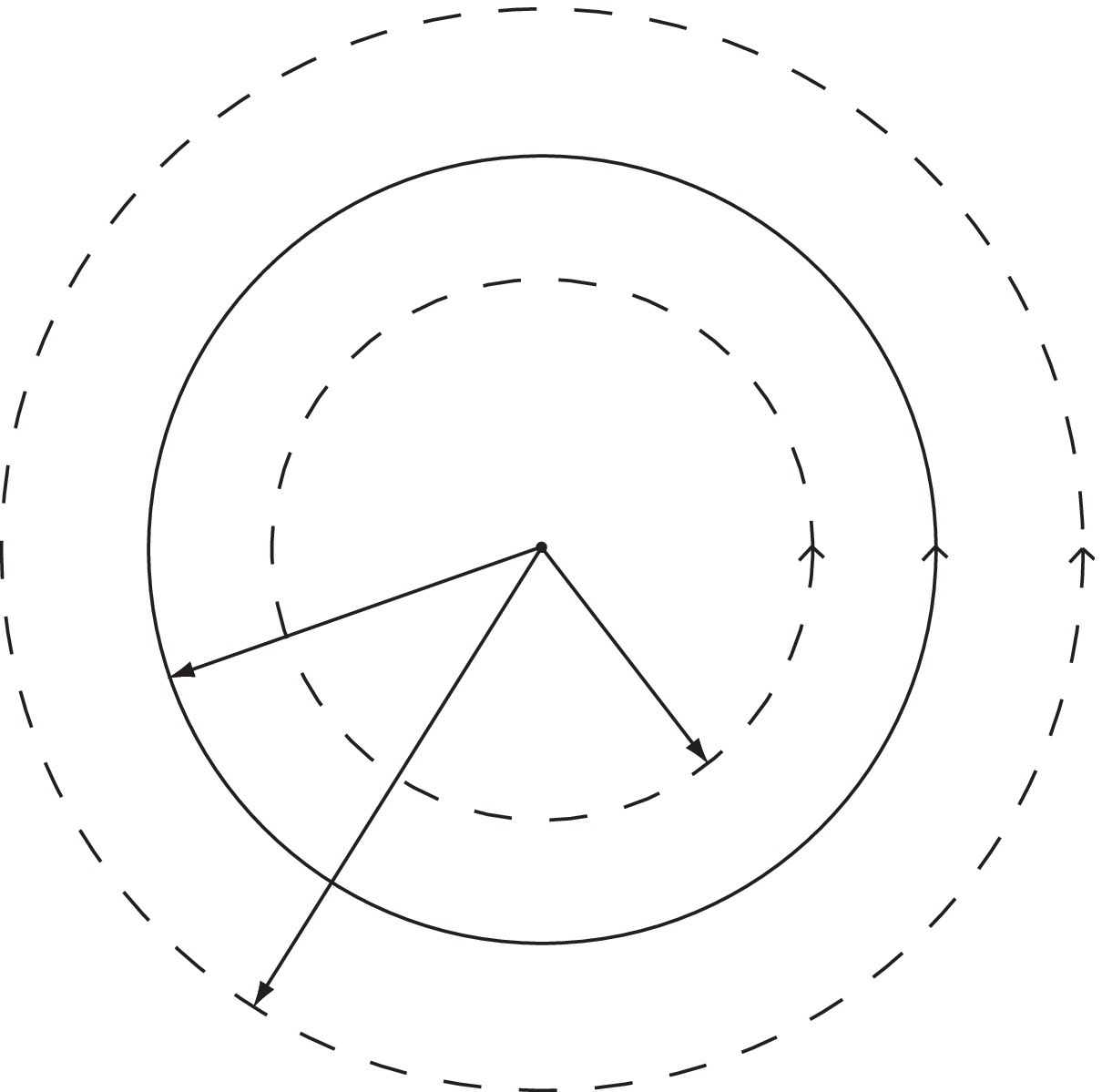}
      \put(79,24){$\T_1$}
         \put(70,32){$\T_r$}
          \put(90,17){$\T_{1/r}$}
          \put(60,60){$\Omega_0$}
          \put(69,69){$\Omega_+$}
          \put(78,78){$\Omega_-$}
          \put(90,90){$\Omega_\infty$}
          \put(54,38){\small $r$}
          \put(27,10){\small $1/r$}
          \put(20,41){\small $1$}
\end{overpic}
\caption{Opening lenses.}\label{fig:nice_case1}
\end{figure}

Let us define
\begin{equation}\label{defH}
H(z)\isdef \begin{cases} \begin{pmatrix} z^{-n} & 0 \\ 0 & z^n
\end{pmatrix}, & \text{if } |z|>1, \\
I, & \text{if } |z|<1,
\end{cases}
\end{equation}
and put
\begin{equation}\label{YtoT}
    T(z) \isdef Y(z)\, H(z)
\,.
\end{equation}
Then $T$ is holomorphic in $\C\setminus {\T_1}$,
$$
T_+(t)=T_-(t)\, \begin{pmatrix} t^n & w(t)  \\ 0 & t^{-n}
\end{pmatrix}\,, \quad t\in \T_1\,,
$$
and
$$
\lim_{z\to \infty} T(z)=I\,.
$$
Next transformation is induced by the following factorization of
the jump matrix for $T$:
 \begin{equation*}\label{equ:factorization}
\begin{pmatrix} t^n & w(t)  \\ 0 & t^{-n}
\end{pmatrix}= \begin{pmatrix} 1 & 0  \\ 1/(t^n w(t)) & 1
\end{pmatrix}\, \begin{pmatrix} 0 & w(t)  \\ -1/w(t) & 0
\end{pmatrix}\, \begin{pmatrix} 1 & 0  \\ t^n/ w(t)  & 1
\end{pmatrix}\,.
\end{equation*}
Hence, with $r$, $\rho <r<1$, fixed above we determine the
following regions (Figure \ref{fig:nice_case1}):
\begin{align*}
 \Omega_0 &= \{z:\, |z|<r \}, \qquad \Omega_\infty =\{z:\,
|z|>1/r \} ,\\
\Omega_+ &=\{z:\, r<|z|<1 \}, \qquad \Omega_- = \{z:\, 1<|z|<1/r
\} .
\end{align*}

We define now
\begin{equation}\label{TtoU}
    U(z) \isdef T(z)K(z),
\end{equation}
where
\begin{equation}\label{KforT}
    K(z) \isdef \begin{cases} I, & \text{if } z \in \Omega_0 \cup \Omega_\infty, \\
\begin{pmatrix} 1 & 0  \\ z^n/  w(z)  & 1
\end{pmatrix}^{-1}\,, & \text{if }
z\in \Omega_+, \\ \begin{pmatrix}  1 & 0  \\ 1/(z^n w(z)) & 1
\end{pmatrix}\,, & \text{if }
z\in \Omega_-.
\end{cases}
\end{equation}
Then $U$ is holomorphic in $\C\setminus \left(  {\T_r} \cup {\T_1}
\cup {\T_{1/r}}\right)$, $$ \lim_{z\to \infty} U(z)=I\,,
$$
and
$$
U_+(t)=U_-(t)\, J_U(t), \quad t\in \left( {\T_r}  \cup {\T_1}\cup
{\T_{1/r}}\right),
$$
where
$$
J_U(t)=\begin{cases} \begin{pmatrix} 0 & w(t)  \\ -1/w(t) & 0
\end{pmatrix}, & \text{if } t\in {\T_1}, \\   \begin{pmatrix} 1 & 0  \\ t^n/  w(t)  & 1
\end{pmatrix}, & \text{if } t\in {\T_r}, \\
\begin{pmatrix} 1 & 0  \\ 1/(t^n w(t)) & 1
\end{pmatrix}, & \text{if } t\in {\T_{1/r}}.
\end{cases}
$$

Using the piece-wise analytic matrix
\begin{equation}\label{equ:defN}
    N(z)\isdef \begin{cases}
   \begin{pmatrix} 1/\DD_{\rm e}(z) & 0  \\ 0 & \DD_{\rm e}(z)
\end{pmatrix}=\begin{pmatrix} D_{\rm e}(z)/\tau  & 0  \\ 0 & \tau /D_{\rm e}(z)
\end{pmatrix}, & \text{if } |z|>1\,, \\
  \begin{pmatrix} 0 & \DD_{\rm i}(z)  \\ -1/\DD_{\rm i}(z) & 0
\end{pmatrix}=\begin{pmatrix} 0 & D_{\rm i}(z)/\tau   \\ -\tau /D_{\rm i}(z) & 0
\end{pmatrix}, & \text{if } |z|<1\,.
    \end{cases}
\end{equation}
we make a new transformation, defining
\begin{equation}\label{equ:defS}
S(z) \isdef U(z) N^{-1}(z)\,.
\end{equation}
Matrix $S$ is holomorphic in $\C\setminus ({\T_r} \cup {\T_1}\cup
{\T_{1/r}} )$, $ \lim_{z\to \infty} S(z)=I$,  and
$$
S_+(t)=S_-(t)\, J_S(t), \quad t\in  {\T_r} \cup {\T_1}\cup
{\T_{1/r}},
$$
where
$$
J_S(t)=N_- J_U N_+^{-1}= \begin{cases} I, & \text{if } t\in {\T_1}, \\
\begin{pmatrix} 1 & -t^n \DD_{\rm i}^2(t)/  w(t)  \\ 0  & 1
\end{pmatrix}, & \text{if } t\in {\T_r}, \\
\begin{pmatrix} 1 & 0  \\ \DD_{\rm e}^2(t)/(t^n w(t)) & 1
\end{pmatrix}, & \text{if } t\in {\T_{1/r}}.
\end{cases}
$$
This shows in particular that $S$ is analytic across ${\T_1}$.
Using the definition of $\FF$, the jump matrix function $J_S$ can
be rewritten as
\begin{equation}\label{def_J_S}
J_S(t)=N_- J_U N_+^{-1}= \begin{cases} I, & \text{if } t\in {\T_1}, \\
\begin{pmatrix} 1 & -t^n \FF(t)/  \tau ^2  \\ 0  & 1
\end{pmatrix}, & \text{if } t\in {\T_r}, \\
\begin{pmatrix} 1 & 0  \\ \tau ^2/(t^n \FF(t)) & 1
\end{pmatrix}, & \text{if } t\in {\T_{1/r}}.
\end{cases}
\end{equation}
Observe that
\begin{equation}\label{boundsForJ_S}
J_S(z)=I+\mathcal O (r^n)\,, \quad z \in \T_r \cup \T_{1/r}\,.
\end{equation}

Unravelling these transformations we have
\begin{equation}\label{equ:exprForY}
    Y(z)=S(z)N(z)K^{-1}(z)H(z)\,,
\end{equation}
where $H$, $K$ and $N$ are defined in \eqref{defH}, \eqref{KforT}
and \eqref{equ:defN}, respectively.

\begin{proof}[Proof of Theorem \ref{thm:nice case}]
In order to make formula \eqref{equ:exprForY} useful we need an
asymptotic expression for $S$. It can be obtained if we observe
that $S$ solves the following Riemann-Hilbert problem:
$S_+=S_-J_S$ on ${\T_r} \cup {\T_{1/r}}$, and $S(z)=I+\mathcal O
(1/z)$ as $z\to \infty$. The jump condition may be rewritten
equivalently as
\begin{equation}\label{RHforS}
(S-I)_+ -(S-I)_-=S_- (J_S-I)\quad \text{on }{\T_r} \cup {\T_{1/r}}
\,.
\end{equation}
Let us introduce the following notation: for $f$ defined and
continuous on $\T_r \cup \T_{1/r}$,
$$
\CC(f)(z) \isdef \frac{1}{2\pi i}\,\oint_{{\T_r} \cup {\T_{1/r}} }
\frac{f(t)\, dt}{t-z}\,,
$$
is the Cauchy operator, and
$$
\CC_-(f)(z) \isdef \lim_{\stackrel{u\to z}{u\in  \text{``$-$''
side}}}\CC(f)(u)\,, \qquad \CC_+(f)(z) \isdef \lim_{\stackrel{u\to
z}{u\in  \text{``$+$'' side}}}\CC(f)(u)\,,
$$
are the Cauchy projection operators on the ``$-$'' and
``$+$''-sides of the corresponding curves, respectively.
Furthermore, let
$$
\PP(f)(z) \isdef \CC(f(J_S-I))(z), \qquad \PP_-(f)(z) \isdef
\CC_-(f(J_S-I))(z)\,.
$$

We can use the fact that Sokhotsky--Plemelj's formula
$\CC_+-\CC_-=I$ is valid also in $L^2$ (cf.\ \cite{Stein70} and
\cite[\S 7.1]{MR2000g:47048}).

Then, due to the normalization at infinity, \eqref{RHforS} can be
expressed, equivalently, as
\begin{equation}\label{S}
S -I=\CC (S_-  (J_S -I)) =\PP(S_-) = \PP(I) +\PP(S_- -I) \,.
\end{equation}
Considering the ``$-$'' boundary values in \eqref{S}, we get that
$\mu \isdef S_- -I$ satisfies $ \mu =\PP_-(I)+\PP_-(\mu) $, or
\begin{equation}\label{equation for mu}
(I-\PP_-)(\mu)=\PP_-(I)\,.
\end{equation}
This may be viewed as a singular integral equation for $\mu$.
Since $\CC_-$ is a bounded linear operator in $L^2({\T_r} \cup
{\T_{1/r}})$, there exists a constant $0<A<+\infty$, such that $\|
\PP_- f\|\leq A \|f\|$ for every $f\in L^2({\T_r} \cup {\T_{1/r}}
)$, where we use the $L^2$ norm. In particular, using
\eqref{boundsForJ_S} we get
$$
\|\PP_-(f)\|\leq A \| f(J_S-I)\|\leq A \|J_S-I\|_{L^\infty} \,
\|f\|\leq A r^n \|f\|\,,
$$
for $n$ sufficiently large. For these $n$'s, $\|\PP_-\|<1$, and we
can solve \eqref{equation for mu} using Neumann series:
$$
\mu=(I-\PP_-)^{-1}(\PP_-(I))=\sum_{k=0}^\infty
\PP_-^k(\PP_-(I))=\sum_{k=1}^\infty \PP_-^k(I)\,,
$$
or equivalently,
$$
S_-=\sum_{k=0}^\infty \PP_-^k(I)\,.
$$
Properly speaking, if we use the $L^{2}$ theory of
Riemann--Hilbert problems, this series can be shown to converge in
the $L^{2}$ sense. However, the series is, in fact, uniformly
convergent.  Indeed, straightforward contour deformations show
that each successive integral in ${\cal P}_{-}^{k}(I)$ taken over
a circle of radius $r$ may be deformed to a circle of radius $r'$
with $r'<r$, and similarly each successive integral taken over a
circle of radius $1/r$ may be deformed to a circle of radius
$1/r'$ with $r'<r$. Then arguments quite analogous to those
leading to \eqref{estimatesFk} may be used to estimate each term
${\cal P}_{-}^{k}(I)$, and the uniform convergence is evident.

Replacing this series in \eqref{S} we get the full asymptotic
expansion for $S$:
\begin{align*}
S(z)& =I +\PP(S_-)= I +\PP\left(\sum_{k=0}^\infty \PP_-^k(I)
\right)= I +\sum_{k=0}^\infty \PP\left( \PP_-^k(I) \right)=
\sum_{k=0}^\infty S^{(k)}\,,
\end{align*}
where
$$
S^{(0)}=I\,, \quad \text{and} \quad S^{(k+1)}=\PP \left(S^{(k)}_-
\right)\,, \quad k\geq 0\,.
$$
Using the explicit expression for $J_S$ we get
$$
S^{(1)}(z)=\frac{1}{2\pi i}\, \oint_{{\T_r} \cup {\T_{1/r}} }
\frac{
(J_S(t)-I)\, dt}{t-z}=\begin{pmatrix} 0 & f_n^{(1)}(z)  \\
g_n^{(1)}(z)  & 0
\end{pmatrix}\,,
$$
where
$$
f_n^{(1)} \isdef \mathcal M_n^{\rm i}(1) \,, \qquad g_n^{(1)}
\isdef \mathcal M_n^{\rm e}(1) \,,
$$
and operators $\mathcal M_n^{\rm i}$, $\mathcal M_n^{\rm e}$ have
been defined in \eqref{MPlus}--\eqref{MMinus}.  Furthermore,
\begin{align*}
S^{(2)}(z)& =\PP(S_-^{(1)})(z)=\frac{1}{2\pi i}\,  \oint_{{\T_r}
\cup {\T_{1/r}} } \frac{
S_-^{(1)}(t)(J_S(t)-I)\, dt}{t-z} \\
 &=\frac{1}{2\pi i}\,
\oint_{{\T_{1/r}} } \begin{pmatrix} 0 & (f_n^{(1)})_-(t)  \\
(g_n^{(1)})_-(t)  & 0
\end{pmatrix}\, \begin{pmatrix} 0 & 0  \\ \DD_{\rm e}^2(t)/(t^n
w(t)) & 0 \end{pmatrix}\frac{  dt}{t-z} \\
 & + \frac{1}{2\pi i}\,
\oint_{{\T_r} } \begin{pmatrix} 0 & (f_n^{(1)})_-(t)  \\
(g_n^{(1)})_-(t)  & 0
\end{pmatrix}\, \begin{pmatrix} 0 & -t^n \DD_{\rm i}^2(t)/  w(t)  \\ 0 & 0 \end{pmatrix}\frac{  dt}{t-z}\\
 &=\frac{\tau ^2}{2\pi i}\,
\oint_{{\T_{1/r}} } \begin{pmatrix} 1 & 0  \\
0 & 0
\end{pmatrix}\,\frac{(f_n^{(1)})_-(t)}{\FF(t) t^n}\, \frac{  dt}{t-z}  - \frac{1}{2\pi i\, \tau ^2}\,
\oint_{{\T_r} } \begin{pmatrix} 0 & 0  \\
0 & 1
\end{pmatrix}\, t^n \FF(t) (g_n^{(1)})_-(t)\,\frac{  dt}{t-z}\\
  &= \begin{pmatrix} f_n^{(2)}(z) & 0  \\
0 & g_n^{(2)}(z)   \end{pmatrix}\,,
\end{align*}
where
$$
g_n^{(2)} \isdef \left(\mathcal M_n^{\rm i} \circ \mathcal
M_n^{\rm e}\right)(1)=  M_n^{\rm i} (g_n^{(1)}) \,, \qquad
f_n^{(2)} \isdef \left(\mathcal M_n^{\rm e} \circ \mathcal
M_n^{\rm i}\right)(1)=  M_n^{\rm e} (f_n^{(1)}) \,.
$$

Iterating this process, we get that
\begin{equation}\label{asymp_fla_for_S}
S =\begin{pmatrix} \displaystyle  1+\sum_{k=1}^\infty f_n^{(2k )}  &
 \displaystyle  \sum_{k=0}^\infty f_n^{(2k+1)}   \\
 \displaystyle  \sum_{k=0}^{\infty} g_n^{(2k+1)}   &
  \displaystyle  1+ \sum_{k=1}^\infty
g_n^{(2k)}
\end{pmatrix},
\end{equation}
where the matrix entries are defined recursively by
$$
f_n^{(1)} \isdef  \mathcal M_n^{\rm i}(1) \,, \quad f_n^{(2)}
\isdef \mathcal M_n^{\rm e} (f_n^{(1)})\,, \quad f_n^{(2k+1)}
\isdef \mathcal M_n^{\rm i} (f_n^{(2k)}) \,, \quad f_n^{(2k+2)}
\isdef \mathcal M_n^{\rm e} (f_n^{(2k+1)})\,, \quad
 k \in
\N\,,
$$
and
$$
g_n^{(1)} \isdef  \mathcal M_n^{\rm e}(1) \,, \quad g_n^{(2)}
\isdef \mathcal M_n^{\rm i}(g_n^{(1)})\,, \quad g_n^{(2k+1)}
\isdef \mathcal M_n^{\rm e} (g_n^{(2k)}) \,, \quad g_n^{(2k+2)}
\isdef \mathcal M_n^{\rm i} (g_n^{(2k+1)})\,, \quad
 k \in
\N\,.
$$

It is usually the case in the asymptotic analysis of
Riemann--Hilbert problems that the expansion for the error matrix
is, at the end of the day, an asymptotic expansion, rather than
being uniformly convergent. However, in the present setting,
estimates \eqref{estimatesFk} and \eqref{estimatesGk} show that
the series appearing in \eqref{asymp_fla_for_S} are actually
uniformly convergent (at least, for $n$ satisfying
\eqref{conditionConvergence})!

Now we may replace \eqref{asymp_fla_for_S} in \eqref{equ:exprForY}
in order to find an expression for $Y$. We must do it
independently in each region.

In the domain  $\Omega_\infty$,
$$ N(z)=   \begin{pmatrix} D_{\rm e}(z)/\tau  & 0  \\ 0 & \tau /D_{\rm e}(z)
\end{pmatrix}, \quad K(z)=I, \quad H(z)=\begin{pmatrix} z^{-n} & 0 \\ 0 & z^n
\end{pmatrix}\,.
$$
Hence,
$$
  Y(z)=S(z)\begin{pmatrix} z^n D_{\rm e}(z)/\tau & 0  \\ 0 & \tau
  /(z^n D_{\rm e}(z))
\end{pmatrix},
$$
and
\begin{equation}\label{asympt_region_omega_inf}
\Phi_n(z)=Y_{11}(z)=S_{11}(z)\,z^n D_{\rm e}(z)/\tau =
 \frac{z^n D_{\rm e}(z)}{\tau}\left(  1+\sum_{k=1}^\infty f_n^{(2k)}(z)\right)\,,
\quad |z|> 1/r\,.
\end{equation}

In the domain $\Omega_-$,
$$ N(z)=   \begin{pmatrix} D_{\rm e}(z)/\tau  & 0  \\ 0 & \tau /D_{\rm e}(z)
\end{pmatrix}, \quad K(z)=\begin{pmatrix} 1 & 0  \\ 1/(z^n w(z)) & 1
\end{pmatrix}, \quad H(z)=\begin{pmatrix} z^{-n} & 0 \\ 0 & z^n
\end{pmatrix}\,.
$$
Hence,
$$
  Y(z)=S(z)\begin{pmatrix} z^n D_{\rm e}(z)/\tau  & 0  \\ -\tau /D_{\rm i}(z) & \tau /(z^n D_{\rm e}(z)
  )
\end{pmatrix},
$$
and for $1\leq |z|< 1/r$,
\begin{equation}\label{asympt_region_omega_minus}
\begin{split}
 \Phi_n(z)&=Y_{11}(z)=S_{11}(z)\,\frac{z^n D_{\rm e}(z)}{\tau }-
S_{12}(z)\,\frac{\tau }{D_{\rm i}(z)} \\ & =\frac{z^n D_{\rm
e}(z)}{\tau } \left(   1+\sum_{k=1}^\infty f_n^{(2k)}(z) \right) -
\frac{\tau }{D_{\rm i}(z)} \, \left( \sum_{k=0}^\infty
f_n^{(2k+1)}(z) \right) \,.
\end{split}
\end{equation}

In the domain  $\Omega_+$,
$$ N(z)=   \begin{pmatrix} 0 & D_{\rm i}(z)/\tau   \\ -\tau /D_{\rm i}(z) & 0
\end{pmatrix}, \quad K^{-1}(z)=\begin{pmatrix}  1 & 0  \\ z^n/  w(z)  & 1
\end{pmatrix}, \quad H(z)=I\,.
$$
Hence,
$$
  Y(z)=S(z)\begin{pmatrix} z^n  D_{\rm e}(z)/\tau   & D_{\rm i}(z)/\tau   \\ -\tau /D_{\rm i}(z)  & 0
\end{pmatrix},
$$
and we recover formula \eqref{asympt_region_omega_minus}. Hence,
this asymptotic representation is still valid for $r< |z|\leq 1$.

Finally, in the domain $\Omega_0$,
$$
N(z)=   \begin{pmatrix} 0 & D_{\rm i}(z)/\tau   \\ -\tau /D_{\rm
i}(z) & 0
\end{pmatrix}, \quad K(z)=I, \quad H(z)=I\,.
$$
Hence,
\begin{equation}\label{Yinside}
  Y(z)=S(z)\begin{pmatrix} 0 & D_{\rm i}(z)/\tau   \\ -\tau /D_{\rm i}(z) & 0
\end{pmatrix},
\end{equation}
and
\begin{equation}\label{asympPhiNicecase}
\Phi_n(z)=Y_{11}(z)=-S_{12}(z)\, \frac{\tau }{D_{\rm
i}(z)}=-\frac{\tau }{D_{\rm i}(z)}\left( \sum_{k=0}^\infty
f_n^{(2k+1)} (z) \right)\,, \quad |z|< r\,.
\end{equation}
This concludes the proof of Theorem \ref{thm:nice case}.
\end{proof}

\begin{proof}[Proof of Theorem \ref{prop:leading_coeff}]
The asymptotic expansion for the leading coefficient $\kappa_n$ is
a straightforward consequence of formula \eqref{kappa} and the
Riemann-Hilbert steepest descent analysis performed above. Indeed,
by \eqref{asymp_fla_for_S} and  \eqref{Yinside},
$$
Y_{21}(z)=- \frac{\tau\, S_{22}(z)}{D_{\rm i}(z)}=- \frac{\tau\,
}{D_{\rm i}(z)}\, \left( 1+ \sum_{k=1}^\infty
g_n^{(2k)}(z)\right)\,,
$$
and the statement follows from \eqref{tau} and \eqref{kappa}.
\end{proof}


\begin{proof}[Proof of Proposition \ref{prop:residues}]
Assume first that $|z|\leq (\rho '+\rho)/2$. Then $t^n \FF(t)/(t-z)$
is a meromorphic function of $t$ in the annulus $(\rho '+\rho)/2
<|t|<r$ whose only poles are $a_1, \dots , a_u$. Hence, taking $s=(3
\rho +\rho ')/4$, we may deform the path of integration used in the
definition of $f_n^{(1)}$ into the one depicted in Figure
\ref{fig:deforming}, left, obtaining that
\begin{figure}[htb]
\centering \begin{overpic}[scale=0.53]{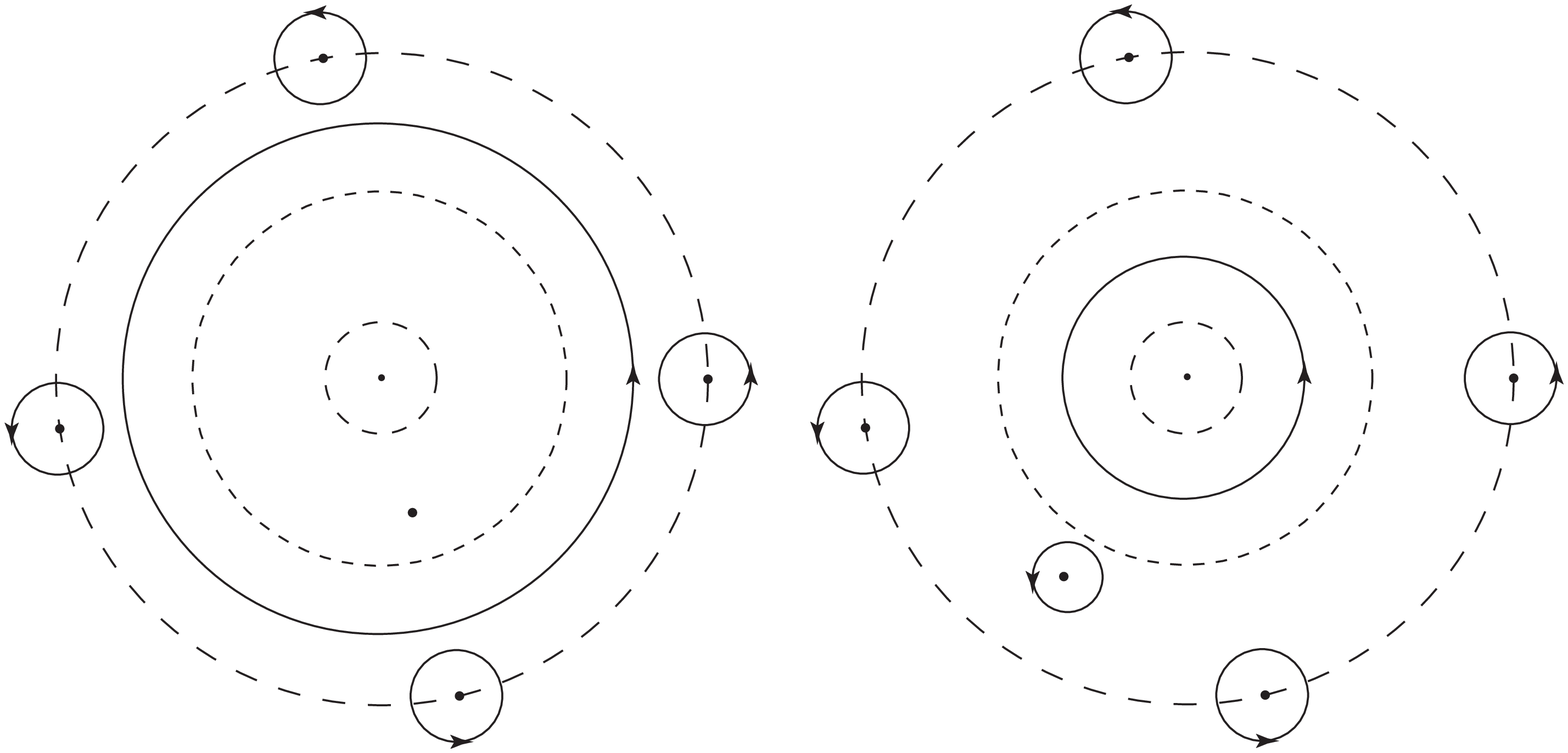}
         \put(30,39.5){\small $\T_{s }$}
         \put(80,30.8){\small $\T_{s}$}
          \put(93,9){\small $\T_{\rho }$}
          \put(41.5,9){\small $\T_{\rho }$}
          \put(27.5,20){\small $\T_{\rho' }$}
         \put(79,20){\small $\T_{\rho' }$}
          \put(42.5,23){\small $a_k$}
          \put(94.2,23){\small $a_k$}
          \put(25,23){\small $0$}
          \put(77,23){\small $0$}
          \put(24.5,14){\small $z$}
          \put(68.5,10){\small $z$}
\end{overpic}
\caption{Case of a finite number of singularities on $\T_\rho
$.}\label{fig:deforming}
\end{figure}
$$
f_n^{(1)}(z)=-\frac{1}{2\pi i\, \tau^2}\, \oint_{{\T_r}}
 \frac{\FF(t) \, t^n }{ t-z} \, dt=-\frac{1}{ \tau^2}\,\sum_{k=1 }^u \res_{t=a_k } \left( \FF (t) \,
\frac{t^n }{ t-z } \right)
   -\frac{1}{2\pi i\, \tau^2}\, \oint_{{\T_{s}}}
 \frac{\FF(t) \, t^n }{ t-z} \, dt\,.
$$
Since $|z|\leq (\rho '+\rho)/2$, we have
\begin{equation}\label{estimateLastTerm}
\left| \oint_{{\T_{\widehat r }}}
 \frac{t^n  \FF(t)}{ t-z} \, dt\right| \leq 4\pi \, \max_{|t|=\widehat r} \FF(t)\, \widehat r^n\,.
\end{equation}
Hence, taking into account \eqref{roughEstimatePhi}--\eqref{boundH}
we obtain \eqref{asymptotics_finite_sing} for $|z|\leq (\rho '+\rho
)/2$, with $\delta =s/\rho $. Obviously, as long as $|z|\leq s$,
\eqref{case_finite_number_sing2} is a consequence of
\eqref{asymptotics_finite_sing}.

Assume now that $(\rho '+\rho)/2<|z|\leq r'  $. Then $t^n
\FF(t)/(t-z)$ is a meromorphic function of $t$ in the annulus
$\rho '<|t|<r$ whose poles are $a_1, \dots , a_u$ and $z$. If $z
\notin \{a_1, \dots, a_u \}$, the pole of $t^n \FF(t)/(t-z)$ at
$t=z$ is simple, and
$$
\res_{t=z } \left( \FF (t) \, \frac{t^n }{ t-z } \right)=z^n \FF
(z)\,.
$$
Deforming the path of integration $\T_r$ for $f_n^{(1)}$ into the
one depicted in Figure \ref{fig:deforming}, right, with $s =(3 \rho'
+\rho )/4$, we obtain that
$$
f_n^{(1)}(z)=-\frac{1}{ \tau^2}\,z^n \FF (z) -\frac{1}{
\tau^2}\,\sum_{k=1 }^u \res_{t=a_k } \left( \FF (t) \, \frac{t^n
}{ t-z } \right)
   -\frac{1}{2\pi i\, \tau^2}\, \oint_{{\T_{s}}}
 \frac{\FF(t) \, t^n }{ t-z} \, dt\,.
$$
However, since $|z|>(\rho +\rho ')/2$, the same estimate
\eqref{estimateLastTerm} for the last term in this identity
remains valid. This proves \eqref{case_finite_number_sing2}, which
remains valid after possible cancellations even for $z \in \{a_1,
\dots, a_u \}$. Finally, for any compact set $K\subset \D_\rho $,
we have that $\max_{z\in K} |z/\rho|\leq \delta<1 $, and
\eqref{asymptotics_finite_sing} follows from
\eqref{case_finite_number_sing2}.
\end{proof}

\begin{proof}[Proof of Theorem \ref{prop:many dominant}]
Observe first that if $a$, $0<|a|<1$, is a pole of an analytic
function $ f$ in $|z|<1$, and in a neighborhood of $z=a$ the
Laurent expansion of $f$ is
$$
f(z)=\frac{f_{-m}}{(z-a)^m}+\frac{f_{-m+1}}{(z-a)^{m-1}}+\dots\,,
\quad f_{-m}\neq 0\,, \quad m \in \N\,,
$$
then
%
\begin{align*}
 f(t) t^n   & =f(t) (t-a+a)^n =f(t)\sum_{j=0}^{n}
\binom{n}{j}\,  a^{n-j} (t-a)^j \,,
\end{align*}
and for $n\geq m$ we get
\begin{equation}\label{formulaforres}
\res_{t=a } \left( f(t) t^n   \right)   = \sum_{j=0}^{m-1}
\binom{n}{j}\,  a^{n-j}  f_{-j-1}= f_{-m}   \, a^{n-m+1}
\,\binom{n }{ m-1 }\, \left(1 + h_n\right) \,,
\end{equation}
where $h_n=0$ for $m=1$, and for $m\geq 2$,
$$
|h_n|\leq \frac{C_1}{n}\, \frac{\max_{0\leq j \leq m-2}
|f_{-j-1}|}{|f_{-m}|}
$$
and constant $C_1$ depends on $a$ and $m$ only.
%

Thus, with assumptions of Theorem \ref{prop:many dominant}, for
$k=1, \dots, \ell$,
\begin{equation}\label{residueDominant}
\res_{t=a_k } \left( \FF (t)   \, \frac{t^n  }{ t-z } \right)=
f_{-m,k}(z)   \, a_k^{n-m+1} \,\binom{n }{ m-1 }\, \left(1 +
h_{n,k}(z)\right)\,,
\end{equation}
with
\begin{equation}\label{residueFF}
f_{-m,k}(z)=\lim_{z\to a_k} \frac{\FF (t)  (t-a_k)^m }{ t-z
}=\frac{D_{\rm i} (a_k)\,  \widehat D_{\rm e} (a_k) }{  a_k-z
 }\,.
\end{equation}
Furthermore, for  $z \notin B_\varepsilon (a_k)$  we have that
$$
|h_{n,k}(z)|\leq \frac{C_2}{n}\, \frac{1}{|z-a_k|^{m-1}}\leq
\frac{C_2}{n}\, \frac{1}{\varepsilon ^{m-1}} \,,
$$
and constant $C_2$ depends on $\rho $, $r$, $\Lambda$ and $m$, but
neither on $z$ nor $n$.

The multiplicity $m'$ of a non-dominant pole $a_k$, $k>\ell$, is
strictly less than $m$; by \eqref{formulaforres} we obtain again
that as long as $z \notin B_\varepsilon (a_k)$,
\begin{equation}\label{residueNonDominant}
\res_{t=a_k } \left( \FF (t)   \, \frac{t^n  }{ t-z } \right)=
\mathcal O \left( \rho ^{n}  n ^{ m'-1 }\right)\,, \quad m'<m\,,
\end{equation}
uniformly in $z$.

Assumptions of Theorem \ref{prop:many dominant} yield those of
Proposition \ref{prop:residues}. Hence, plugging formulas
\eqref{residueDominant}--\eqref{residueNonDominant} into
\eqref{case_finite_number_sing2}--\eqref{asymptotics_finite_sing},
we obtain
\eqref{case_finite_number_dominant_sing}--\eqref{newBOundForH}.
Estimate \eqref{asymptotics_dominant_poles} is a straightforward
consequence of \eqref{case_finite_number_dominant_sing} since we
may assume now $r^3<\rho <r$.

For the last statement it is sufficient to observe that given a
compact subset $K$ in $|z|<\rho $, by
\eqref{asymptotics_dominant_poles},
$$
\frac{1}{a_1 ^{n-m+1}}\, \binom{n }{ m-1 }^{-1}\,\Phi_n(z)
$$
is a normal family in $|z|<\rho $, and every limit of a convergent
subsequence is a rational function with at most $\ell-1$ zeros on
$K$.
\end{proof}

\begin{proof}[Proof of Theorem \ref{thm:clock}]
We carry out the proof in the spirit of \cite{Simon04b}, but with
an additional difficulty generated by the ``floating'' zeros of
$G_n$.

Let  $\varepsilon >0$; without loss of generality we may assume
that $a_1=\rho>0 $. Then for $z\in \mathcal B(\varepsilon )$ we
may rewrite \eqref{case_finite_number_dominant_sing} using the
notation \eqref{G_n}:
\begin{equation}\label{case_finite_number_dominant_singRewritten}
\frac{\tau }{D_{\rm e}(z) }\, \rho
^{-n}\,\binom{n}{m-1}^{-1}\,\Phi_n(z)=\binom{n}{m-1}^{-1}\, \left(
\frac{z}{\rho }\right)^n +\frac{G_n(z)}{\rho ^{m-1}\, \FF(z)}   +
H_n(z) \,,
\end{equation}
where
$$
 H_n(z) =\frac{\tau }{D_{\rm e}(z) }\, \rho
^{-n}\,\binom{n}{m-1}^{-1}\,h_n(z)\,.
$$
According to \eqref{newBOundForH}, there exists a constant
$C=C(\varepsilon )$ such that
\begin{equation}\label{boundForH}
|H_n(z)|\leq \frac{C}{n}\,, \quad  z\in \mathcal B(\varepsilon
)\,.
\end{equation}
Furthermore, $\{ G_n\}$ is an equicontinuous family of rational
functions of a bounded degree. In consequence, there exist
constants $q=q(\varepsilon )$ and $Q=Q(\varepsilon )$ such that
for all $n\in \N$,
$$
z\in \mathcal B(\varepsilon )\setminus \mathcal Z_\varepsilon(G_n
) \quad \Rightarrow \quad 0<q \leq  \left|\frac{G_n(z)}{\rho
^{m-1}\, \FF(z)}\right|\leq Q<+\infty\,,
$$
and a constant $\delta =\delta (\varepsilon )>0$ such that
\begin{equation}\label{smallBoundary}
z, \zeta \in \mathcal B(\varepsilon )\setminus \mathcal
Z_\varepsilon(G_n ) , \; |z-\zeta |<\delta  \quad \Rightarrow
\quad \left|\frac{G_n(z)}{\rho ^{m-1}\, \FF(z)}-\frac{G_n(\zeta
)}{\rho ^{m-1}\, \FF(\zeta )}\right|\leq \frac{q}{4}\,.
\end{equation}
For $v_n$ defined in \eqref{defVn}, there exists $n_0\in\N$ such
that for $n\geq n_0$,
\begin{equation}\label{selectionEpsilon}
 v_n\leq r-\varepsilon  , \quad \frac{C}{n}\leq \frac{q}{4}\,, \quad \text{and}
\quad  \left( 1-\delta /2\right)^n<\frac{q}{2} <\frac{3Q}{2}<
\left( 1+\delta /2\right)^n \,.
\end{equation}

Fix $n\geq n_0$ and assume that $\zeta$, $|\zeta |=v_n $, is such
that $B_{ \delta }(\zeta )\subset \mathcal B(\varepsilon )\setminus
\mathcal Z_\varepsilon(G_n ) $. Denote
$$
\frac{G_n(\zeta )}{\rho ^{m-1}\, \FF(\zeta )}=s \, e^{i\alpha} \quad
\text{and} \quad \Delta_n(z)\isdef \frac{G_n(z )}{\rho ^{m-1}\,
\FF(z )}-s \, e^{i\alpha}\,.
$$
Choose a neighborhood of $t=\zeta $ of the following form:
\begin{equation*}
\begin{split}
I_n & =I_n(k_1,k_2)  \\ & \isdef \left\{z\in \C:\, v_n(1-\delta/2)
<|z|<v_n(1 + \delta/2) , \; \frac{\alpha +2 k_1 \pi}{n}<\arg\left(z
\right)<\frac{\alpha +2 k_2 \pi}{n} \right\}\,,
\end{split}
\end{equation*}
where integers $k_1 \leq k_2$ satisfy $\zeta\in I _n$ and
$k_2-k_1\leq \delta n/8$. Then $I _n(k_1,k_2)\subset \mathcal
B(\varepsilon )\setminus \mathcal Z_\varepsilon(G_n ) $. By
\eqref{smallBoundary}, $\left|\Delta_n(z) \right|<s/4$ on the
boundary $\partial I _n$. In consequence,
\begin{equation}\label{smallFunction}
z\in \partial I _n \quad \Rightarrow \quad |\Delta_n(z)   +
H_n(z)|\leq |\Delta_n(z)|   + |H_n(z)|<\frac{s}{2}\,,
\end{equation}
where we have taken into account \eqref{boundForH} and
\eqref{selectionEpsilon}.

With the notation introduced above we rewrite
\eqref{case_finite_number_dominant_singRewritten} as
\begin{equation}\label{case_finite_number_dominant_singRewrittenBis}
\frac{\tau }{D_{\rm e}(z) }\, \rho
^{-n}\,\binom{n}{m-1}^{-1}\,\Phi_n(z)=b_n(z)+\Delta_n(z)   +
H_n(z) \,,
\end{equation}
where
$$
b_n(z)=\binom{n}{m-1}^{-1}\, \left( \frac{z}{\rho }\right)^n +s \,
e^{i\alpha }\,.
$$
For $|z|=v_n(1-\delta/2)$,
$$
|b_n(z)|\geq s-   \left( 1-\delta /2\right)^n >s/2\,,
$$
where we have taken into account a condition in
\eqref{selectionEpsilon}. Analogously, when $|z|=v_n(1 + \delta
/2)$,
$$
|b_n(z)|\geq   \left( 1+\delta /2 \right)^n -s>s/2\,,
$$
again by \eqref{selectionEpsilon}. Finally, when $\arg(z)=(\alpha
 +2 k   \pi)/n$, $k\in \Z$,
$$
\arg\left( \binom{n}{m-1}^{-1}\, \left( \frac{z}{\rho
}\right)^n\right)=\alpha  +2k  \pi \quad \Rightarrow \quad
|b_n(z)|=\binom{n}{m-1}^{-1}\, \left( \frac{|z|}{\rho }\right)^n +s
>s\,.
$$
In other words, we have that $|b_n(z)|>s/2$ on $z\in \partial I
_n$. Taking into account \eqref{smallFunction}, by Rouche's
theorem, $b_n$ and
$$
\frac{\tau }{D_{\rm e}(z) }\, \rho
^{-n}\,\binom{n}{m-1}^{-1}\,\Phi_n(z)
$$
have the same number of zeros in $I _n(k_1,k_2)$. But every pie
slice $I _n(k,k+1)$ of angle $2\pi/n$ contains only one zero of
$b_n$, so in general $\Phi_n$ has exactly $k_2-k_1$ zeros in every
pie slice $I _n(k_1,k_2)\subset \mathcal B(\varepsilon )\setminus
\mathcal Z_\varepsilon(G_n )$ of angle $2\pi(k_2-k_1)/n$.

On the other hand, if $z_0$ is a zero of $\Phi_n$, by
\eqref{case_finite_number_dominant_singRewrittenBis},
$b_n(z_0)=-\Delta_n(z_0)   - H_n(z_0)$, or
$$
\binom{n}{m-1}^{-1}\, \left( \frac{z_0}{\rho }\right)^n =-s \,
e^{i\alpha }-\Delta_n(z_0)   - H_n(z_0)
$$
so that
$$
\binom{n}{m-1}^{-1}\, \left( \frac{|z_0|}{\rho }\right)^n  = s
+\mathcal O\left(\frac{1}{n} \right)\quad \Rightarrow \quad
|z_0|=\rho \left(1+\frac{1}{n}\, \log \binom{n}{m-1}+\mathcal
O\left(\frac{1}{n}\right) \right)\,.
$$
Analogously, if $z_0, z_1$ are the zeros of $\Phi_n$ in two
consecutive pie slices, $I_n(k-1,k)$ and $I_n(k,k+1)$, then
$$
\arg(z_1)-\arg(z_0)=\frac{2\pi}{n}+\mathcal O \left( \frac{1}{n^2}
\right)\,.
$$
\end{proof}

\begin{proof}[Proof of Proposition \ref{prop:szabados}]
By \eqref{asymptotics_dominant_poles},
$$
    Z \subset \bigcap_{k\geq 1} \overline{\bigcup_{n\geq k} \mathcal Z( G_n)
    }\,,
$$
with $G_n$ defined in \eqref{G_n}, and it is sufficient to describe
all the possible limit points of $\{ G_n\} $. Observe that with the
notation introduced before,
\begin{align*}
 G_n(z)& = \sum_{k=1 }^\ell \frac{D_{\rm i} (a_k)\, \widehat D_{\rm
e} (a_k ) }{ a_k-z
 }\, \exp\left(2\pi i (n-m+1) \theta_k\right)\\
  & = \sum_{k=1 }^\ell \frac{D_{\rm i} (a_k)\, \widehat D_{\rm e}
(a_k ) }{ a_k-z
 }\, \exp\left(2\pi i  \sum_{j=1}^v r_{kj}\, (n-m+1)
 \theta_j\right)\,.
\end{align*}

The case $v=1$ is trivial; assume $v\geq 2$. By Kronecker's theorem
(also known as Kronecker-Weyl theorem, see e.g.~\cite[Ch.\
III]{Cassels:1957}), since $\theta_1=1, \theta_2, \dots, \theta_v$
are rationally independent then for any real numbers $X_2, \dots, X
_v$ there exists a sub-sequence $\{n_j\} \subset \N$ such that
$$
\lim_{n_j} e^{ 2\pi i \left((n_j-m+1) \theta _k-\beta _k
\right)}=1\,, \quad k=2, \dots, v\,,
$$
and $(n_j-m+1)p_k/q_k$, $k=1, \dots, \ell$, have limits modulo $\Z$.
This shows that the set of limit points of $\{ G_n\} $ is
parameterized by $0\leq s_k<q_k$, $s_k\in \Z$, $k=1,\dots,\ell$ and
$X_2, \dots, X_v\in \R$, and is given in the left hand side of
\eqref{description_Szabados2}. Now the statement follows for $v\geq
2$.
\end{proof}

\section{Some examples} \label{sec:examples}

We illustrate our results with some examples. Estimates
\eqref{roughEstimatePhi}--\eqref{boundH} justify concentrating on
the asymptotics of $f_n^{(1)}$ and of its zeros in $\D_r$.

Consider first the simplest case of a polynomial weight inducing a
single pole of $D_{\rm e}$ inside $\T_1$. Let $w(z)=|z-a|^2$,
$|z|=1$, where without loss of generality we may assume $a>0$,
$a\neq 1$ (the case $a=0$ is trivial). Then the normalized
representation of $w$, according to Definition
\ref{def:normalrepresentation}, is
$$
w(z)=\left| \frac{z-a}{z} \right|^2, \quad \text{if } a<1, \quad
\text{and} \quad w(z)=\left| 1 -a z \right|^2, \quad \text{if } a
> 1\,.
$$
Assume $ 0< a <1$; in this case
$$
D_{\rm i}(z)= 1-a z\,, \quad D_{\rm e}(z)= \frac{z}{z-a}\,, \quad
\tau =1\,, \quad \text{and }\quad \FF(z)=\frac{z\, (1 - a
z)}{z-a}\,.
$$
In consequence, with $a<r<1$,
$$
f_n^{(1)}(z)= \begin{cases}   \dfrac{  z^{n+1} (1 - a\, z)
-a^{n+1}(  1 -a^2 )}{ a-z } , & \text{ if }
|z|<r\,, \\%
a^{n+1}\, \dfrac{    1 -a^2 }{ z-a } , & \text{ if } |z|>r\,.
\end{cases}
$$
Observe that for $|z|<r$ this is a polynomials of degree $n$, such
that $f_n^{(1)}(a)\neq 0$, and whose zeros distribute
asymptotically uniformly on $\T_{|a|}$, leaving a gap at $z=a$.

Furthermore,
$$
f_n^{(2)}(z)=a^{2n} (1-a^2)(1-a z)\,, \quad z\in \D_{1/r}\,,
$$
and in consequence,
$$
f_n^{(3)}(z)= a^{2n} (1-a^2)\,  \dfrac{  z^{n+1} (1 - a\, z)^2
-a^{n+1}( 1 -a^2 )^2}{ a-z } , \quad z\in \D_r\,.
$$
Then by \eqref{representation},
\begin{align*}
\Phi_n(z)= & -\frac{1}{1-a z}\,  \left( \dfrac{  z^{n+1} (1 - a\,
z) -a^{n+1}(  1 -a^2 )}{ a-z } \right. \\ & \left.+ a^{2n}
(1-a^2)\, \dfrac{ z^{n+1} (1 - a\, z)^2 -a^{n+1}( 1 -a^2 )^2}{ a-z
} +\mathcal O (a^{5n})\right) \,,
\end{align*}
which is valid locally uniformly in $\D_{|a|+\varepsilon }$, with
$\varepsilon >0$. Obviously, we can get further terms in this
asymptotic expansion.

Analogously,
$$
g_n^{(1)}(z)=a^n \, \frac{1-a^2}{1-a z}\,, \quad z \in \D_{1/r}\,,
$$
and
$$
g_n^{(2)}(z)=-a^n (1-a^2) \, \frac{z^{n+1}- a^{n+1}}{z-a} \,,
\quad z \in \D_{r}\,.
$$
Hence, using Theorem \ref{prop:leading_coeff}, we have the
following asymptotics for the leading coefficient $\kappa_n$ of
orthonormal polynomials $\varphi_n$:
$$
\kappa_{n}^2=\frac{1}{2 \pi}\, \left(1-a^{2 n} (1-a^2) + \mathcal
O (r^{4 n}) \right)\,, \quad n\to \infty\,.
$$

\medskip

Much richer is the case of a polynomial weight inducing two
dominant poles of $D_{\rm e}$ inside $\T_1$. Consider
$$
w(z)=\left|(z-a)(z-b)\right|^2, \quad 0<a=|b|<1\,.
$$
Then
$$
D_{\rm i}(z)=(1 - a \, z)(1 -\overline{b}  \, z), \quad D_{\rm
e}(z)=\frac{z^2}{(z-a)(z-b)}\,, \quad \tau =1\,,
$$
and
$$
\FF(z)=z^2\frac{(1 - a \, z)(1 -\overline{b}  \,
z)}{(z-a)(z-b)}.
$$
Let us denote $b=a\, e^{i \pi \theta}$, and fix $a<r<1$. Then we
have for $z\in \D_r$,
\begin{equation}\label{f1ForRational}
f_n^{(1)}(z)=-a^n \left\{\FF(z) \, \left(\frac{z}{a} \right)^{n} +
a \, \frac{1-a^2}{1-e^{i \pi \theta}}\, \left[\frac{1-a^2 e^{-i
\pi \theta}}{a-z} - e^{i \pi (n+2) \theta }\, \frac{1-a^2 e^{ i
\pi \theta}}{b-z} \right] \right\}\,.
\end{equation}
Hence, for $|z| < |a| $, the zeros of $f_n^{(1)}$ are
asymptotically close to the solutions of
$$
\frac{1-a^2 e^{-i \pi \theta}}{a-z}=e^{i \pi (n+2) \theta }\,
\frac{1-a^2 e^{ i \pi \theta}}{b-z}\,,
$$
or in other words, close to
\begin{equation}\label{accumulationZeros}
z=a\, e^{ i \pi \theta}\, \frac{e^{ i \pi (n+1) \theta}-A}{e^{ i
\pi (n+2) \theta}-A}\,, \quad \text{where}\quad  A\isdef
\frac{1-a^2 e^{-i \pi \theta}}{1-a^2 e^{ i \pi \theta}}=
\frac{\overline{1-a^2 e^{ i \pi \theta}}}{1-a^2 e^{ i \pi \theta}}
\,, \quad |A|=1\,.
\end{equation}

\begin{figure}[htb]
\centering \begin{overpic}[scale=0.95]{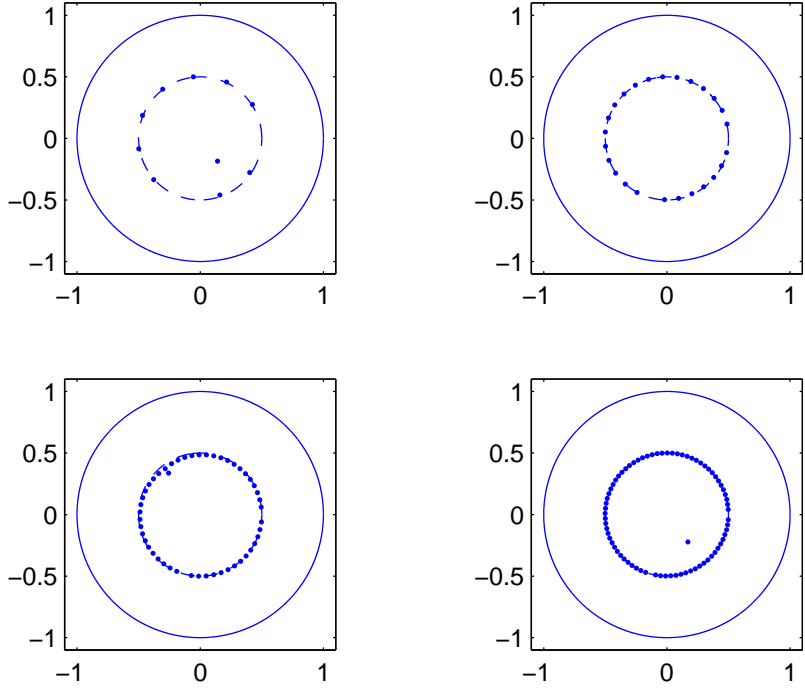}
\end{overpic}
\caption{Zeros of $\Phi_n$ for $n=10, 25, 50, 75$  (from left to
right and from top to bottom) with
$w(z)=|(z-1/2)(z-\zeta_1/2)|^2$, $z\in \T_1$, where
$\zeta_1=\exp(\pi i \sqrt{2})$. The dashed line is the circle of
radius $0.5$.}\label{fig:zerosSimplePoles2}
\end{figure}

Let us consider some particular cases. First,
$\theta=\sqrt{2}\notin \Q$ illustrates Remark \ref{remark:cirves}.
In this case $Z\cap \D_{|a|}$ is parameterized by
$$
z=a\,   \frac{e^{ i \pi \alpha}-A}{ e^{ i \pi \alpha}-e^{ -i \pi
\theta}\, A} \,, \quad \alpha \in \R\,,
$$
which is a straight line passing through the origin  (cf.\ the
numerical example in Fig. \ref{fig:twoandthreepoles}, left). The
rest of the zeros tends to the circle $\T_{1/2}$ exhibiting the
equidistribution pattern predicted by Theorem \ref{thm:clock} (see
Fig. \ref{fig:zerosSimplePoles2}).


Next, if $\theta =1/2$ ($b= ia$), we have by Proposition
\ref{prop:szabados} that the zero accumulation set $Z\cap D_{|a|}$
is finite. Namely, points \eqref{accumulationZeros} in this case
have the following form:
\begin{align*}
 n=4k, \; k\in \N, \quad & \Rightarrow \quad z_0=a\, i\,
 \frac{i-A}{-1-A}\in \D_{|a|} \text{ if } \arg (A) \in (-\pi/4, 3\pi/4) \mod 2\pi\,, \\
  n=4k+1, \; k\in \N, \quad & \Rightarrow \quad z_1=a\, i\, \frac{-1-A}{-i-A}
  \in \D_{|a|} \text{ if } \arg (A) \in ( \pi/4, 5\pi/4) \mod 2\pi\,, \\
  n=4k+2, \; k\in \N, \quad & \Rightarrow \quad z_2=a\, i\, \frac{-i-A}{1-A}
  \in \D_{|a|} \text{ if } \arg (A) \in (3\pi/4, 7\pi/4) \mod 2\pi\,, \\
  n=4k+3, \; k\in \N, \quad & \Rightarrow \quad z_3=a\, i\, \frac{ 1-A}{ i-A}
  \in \D_{|a|} \text{ if } \arg (A) \in (5\pi/4, 9\pi/4) \mod
  2\pi\,.
\end{align*}
\begin{figure}[htb]
\centering \begin{overpic}[scale=0.7]{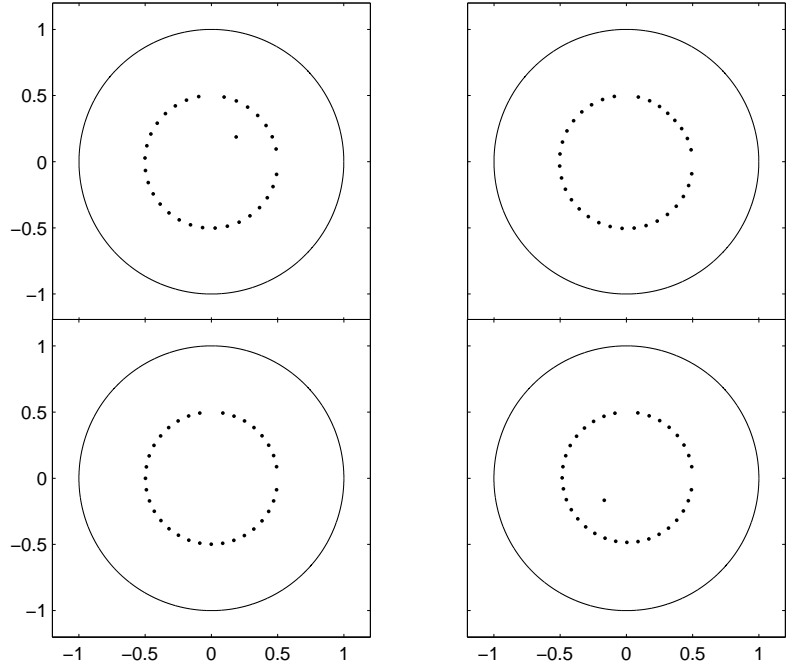}
\end{overpic}
\caption{Zeros of $\Phi_n$ for $w(z)=|(z-1/2)(z-i/2)|^2$, $z\in
\T_1$, and $n=32, \dots, 35$ (from left to right and from top to
bottom).}\label{fig:periodic}
\end{figure}

Hence, $Z\cap \D_{|a|}$ has at most 2 points, and unless $\theta
=\pm \pi/4 \mod \pi$, the sequence $\{f_{n}^{(1)}\}$ will have
zeros in $\D_{|a|}$ for two consecutive indices $n$, and has no
zeros inside $\T_{|a|}$ for the next two. For instance, with
$a=1/2$, $\arg (A)\approx 0.156 \pi$, and zeros of
$\{f_{n}^{(1)}\}$ (and also zeros of $\{\Phi_n \}$) have limits in
$\D_{|a|}$ only for $n=4k$ and $n=4k+3$, $k\in \N$, which are
\begin{equation}\label{spuriousZero}
s_1=\frac{3}{16}\, (1+i)\,, \quad  s_2=-\frac{1}{6}\, (1+i)\,,
\end{equation}
(compare with the numerical experiment in Fig.\
\ref{fig:periodic}).

On the other hand, by \eqref{roughEstimatePhi}--\eqref{boundH} and
by \eqref{f1ForRational},
\begin{equation}\label{limitsSaff}
\frac{\Phi_n(0)}{\Phi_{n-1}(0)} \sim
\frac{f_{n}^{(1)}(0)}{f_{n-1}^{(1)}(0)}=a\, \frac{e^{i \pi (n+2)
\theta}-A}{e^{i \pi (n+1) \theta}-A}\, \left(1+\mathcal O (a^{2n})
\right)\,, \quad n \to \infty\,.
\end{equation}
Let us compare it with the results of
\cite{barrios/lopez/saff:2001} for $\theta=1/2$. Denoting
$$
b_j \isdef \lim_{n=j \mod 4} \frac{\Phi_n(0)}{\Phi_{n-1}(0)}\,,
\quad j=1, \dots, 4\,,
$$
we have by \eqref{limitsSaff},
$$
b_1=a \frac{-i-A}{-1-A}\,, \quad  b_2=a \frac{1-A}{-i-A}\,, \quad
b_3=a \frac{i-A}{1-A}\,, \quad  b_4=a \frac{-1-A}{i-A}\,.
$$
Observe that $|b_1\dots b_4|^{1/4}=|a|$. Using the notation of
\cite{barrios/lopez/saff:2001}, we compute
$$
\Delta_3^{(j)}(z) =\begin{vmatrix} z+b_j & a b_{j+1} & 0 \\
1 & a+ b_{j+1} & z b_{j+2} \\ 0 & 1 & z+ b_{j+2}
\end{vmatrix}\,, \quad j=1, \dots, 4\,,
$$
(where we take $b_5=b_1$ and $b_6=b_2$), and take
$$
\Delta=\bigcup_{j=1}^4 \{z:\, \Delta_3^{(j)}(z) =0\}\,.
$$
For $a=1/2$ it is straightforward to compute that $ \Delta \cap
\D_{|a|}=\{s_1, s_2 \}$, where $s_1$, $s_2$ are defined in
\eqref{spuriousZero}. In other words, this case gives an example
when the whole set $ \Delta \cap \D_{|a|}$ gives the accumulation
points of zeros of $\{ \Phi_n\}$ (cf.\ Remark 1 in
\cite{barrios/lopez/saff:2001}).

%
%

Finally, we dwell very briefly on the case when $D_{\rm e}$ has an
essential singularity on $\T_\rho $ analyzing two examples.

Let first the weight function be given by
\begin{equation}\label{ess_sing}
w(t)=\left|\exp\left( \frac{1}{a-t}\right) \right|^2\,, \quad t \in
\T_1\,,
\end{equation}
with $a\in \D_1$. Without loss of generality we may assume $0<a<1$,
so that $\rho =a$. Then
$$
D_{\rm i}(t)=\exp\left( \frac{t}{a t-1}\right)\,, \quad D_{\rm
e}(t)=\exp\left( \frac{1}{t-a}\right)\,, \quad \tau=D_{\rm
e}(\infty)=1\,,
$$
and
$$
\FF(t)=\exp\left( \frac{t}{ a t-1 }+ \frac{1}{t-a}\right)\,.
$$
Again wee have to analyze the asymptotic behavior of the sequence of
Cauchy transforms
\begin{equation}\label{usePsi}
f_n^{(1)}(z)=  -\frac{1}{2\pi i }\, \oint_{{\T_r}}
 \frac{t^n \FF(t) }{ t-z} \, dt=-\frac{1}{2\pi i }\, \oint _{\T_r} \frac{e^{n \Psi
_n(t)}}{t-z}  \, dt\,, \quad z \in \D_r\,,
\end{equation}
where we have used the notation
\begin{align*}
\Psi _n(t) &\isdef  \log t + \frac{1}{n}\, \log \FF(z)=\log t
+\frac{1}{n}\, \left( \frac{t}{at-1}+\frac{1}{t-a}\right) \,.
\end{align*}
Observe that we cannot apply the classical steepest descent method
to the last integral in \eqref{usePsi}, since $\Psi _n$ depends on
$n$.

Let $\gamma_n$ be a positively oriented closed Jordan curve
encircling $a$ and contained within $\D_r$, and eventually depending
on $n$. Define
\begin{equation}\label{defI}
I_n(z)\isdef -\frac{1}{2\pi i }\, \oint _{\gamma_n } \frac{e^{n \Psi
_n(t)}}{t-z}\, dt\,.
\end{equation}
Then
\begin{equation}\label{f_n}
f_n^{(1)}(z)= \begin{cases} \displaystyle -z^n \FF(z)+I_n(z)\,, & \text{if $z$ is in $\D_r$ and outside $\gamma_n$,} \\
\displaystyle I_n(z)\,, & \text{if $z$ is inside $\gamma_n$.}
\end{cases}
\end{equation}
For all sufficiently large $n$, function $\Psi _n$ has four simple
saddle points, $t_\pm$ and $1/t_\pm$,
$$
0<t_-<a<t_+ <1<1/t_+<1/a<1/t_-\,,
$$
such that
\begin{equation*}
    t_{\pm}= a \pm \frac{\sqrt{a}}{\sqrt{n}}+\mathcal O(1/n)\,, \quad n \to
    \infty\,.
\end{equation*}

\begin{figure}[htb]
\centering \begin{overpic}[scale=0.8]{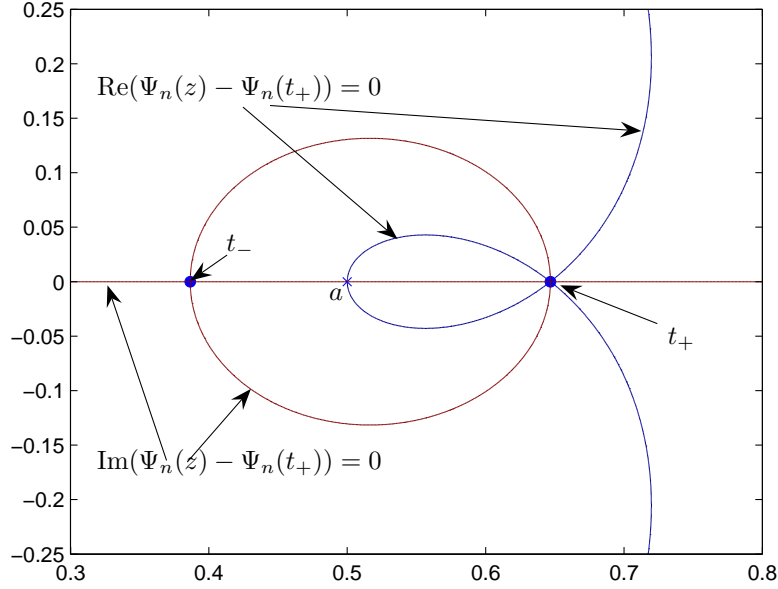}
     \put(16,60){\small $\Re (\Psi_n(z)-\Psi_n(t_+))=0$}
     \put(16,18){\small $\Im (\Psi_n(z)-\Psi_n(t_+))=0$}
     \put(30.5,42.5){\small $t_-$}
     \put(80,32){\small $t_+$}
     \put(42,37){\small $a$}
\end{overpic}
 \caption{Level
curves $\Re (\Psi _n(z)- \Psi _n(t_+))=0$ and $\Im ( \Psi _n(z)-
\Psi _n(t_+))=0$ in a neighborhood of $t_\pm$ (for $a=1/2$ and
$n=30$).}\label{fig:levelsPhi}
\end{figure}

\begin{figure}[htb]
\centering
\begin{tabular}{ll}
\hspace{-1.5cm}\mbox{\begin{overpic}[scale=0.65]{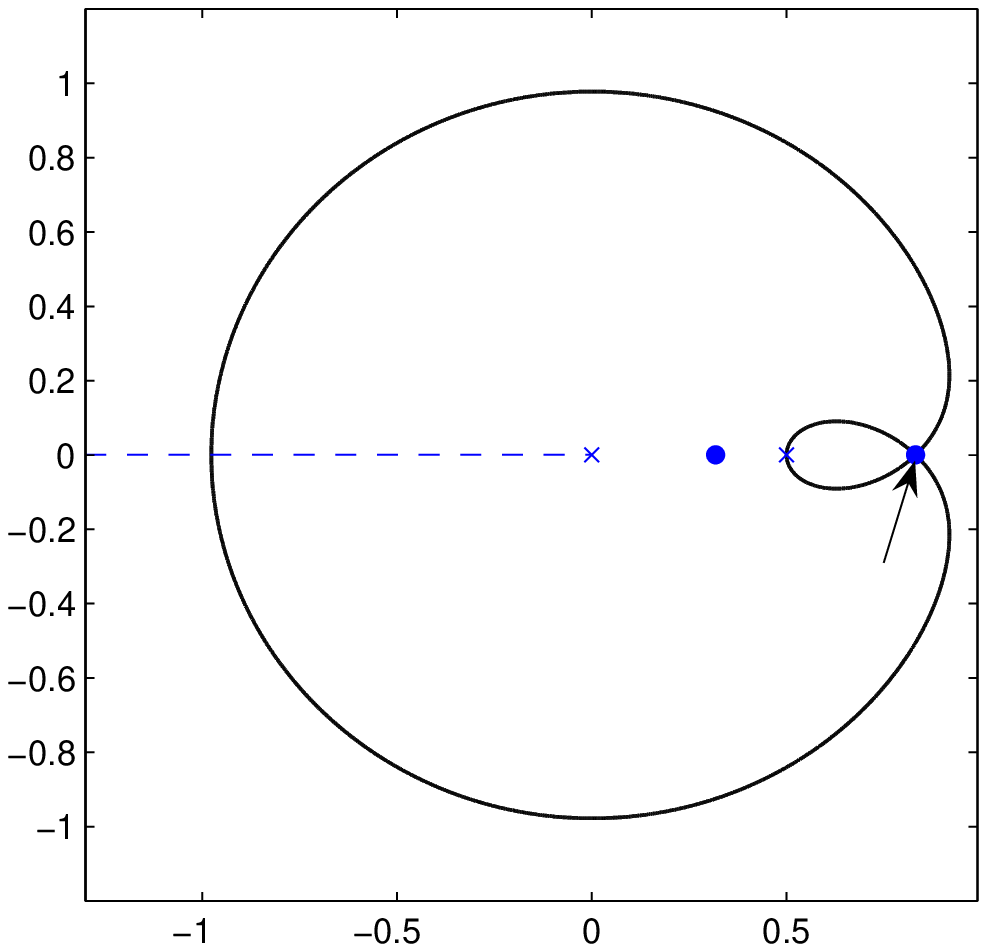}
         \put(54,36.5){\scriptsize $0$}
         \put(70.5,38.5){\scriptsize $a$}
         \put(74,28.5){\scriptsize $t_+$}
\put(36,46){\small $\Re (\Psi_n(z)-\Psi_n(t_+))<0$}
\end{overpic}} &
\hspace{-2cm}\mbox{\begin{overpic}[scale=0.65]{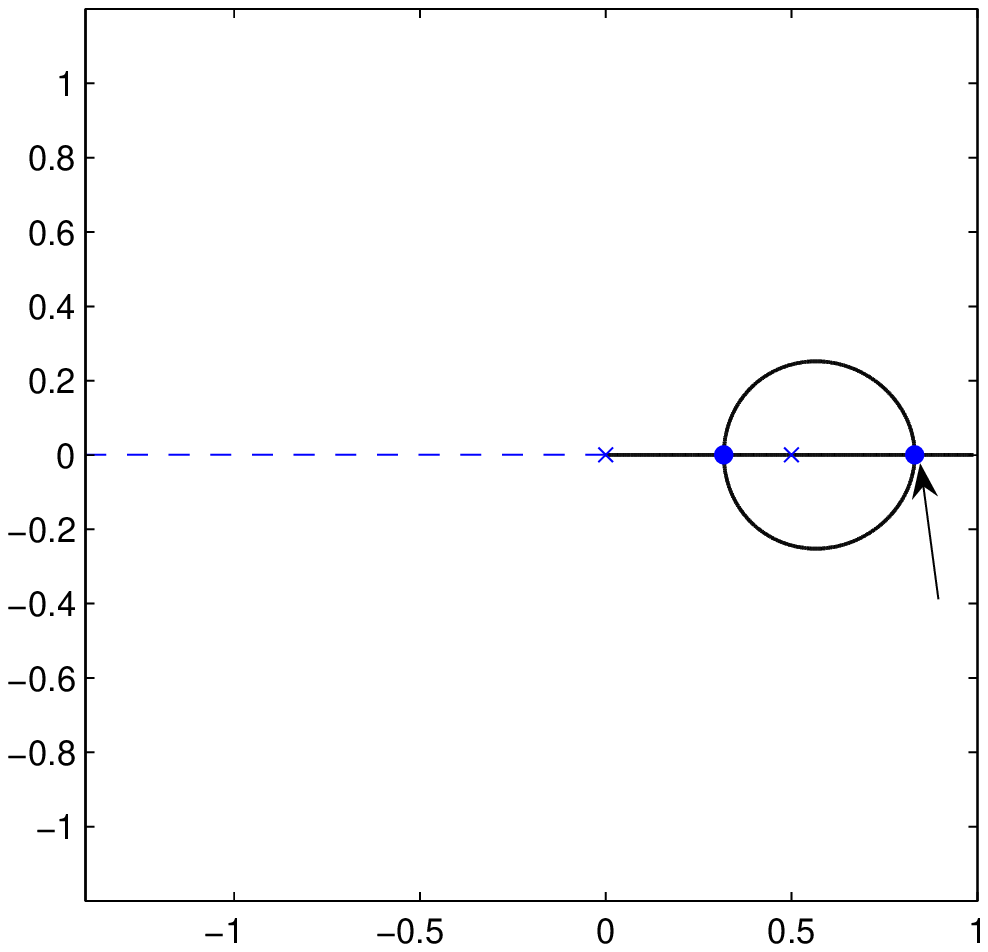}
         \put(54,36.5){\scriptsize $0$}
         \put(69,37){\scriptsize $a$}
         \put(78,26){\scriptsize $t_+$}
\end{overpic}}
\end{tabular}
\caption{Global structure of the level curves $\Re (\Psi _n(z)- \Psi
_n(t_+))=0$ (left) and $\Im (\Psi _n(z)- \Psi _n(t_+))=0$ for
$a=1/2$ and $n=30$.}\label{fig:levelsPhi_dominant_Global}
\end{figure}

Function $\Psi _n$ is meromorphic in the cut disc $\D_1 \setminus
(-1, 0]$, and its level curves $\Gamma=\{t \in \C \setminus
(-\infty, 0]:\, \Im \Psi _n(t)=\Im \Psi _n(t_+)=0\}$ are
trajectories of the quadratic differential $\Psi'_n(t)dt^2$ which
has a double pole at $t=0$ with real leading coefficient, a pole of
order $4$ at $a$, and double zeros at $t_\pm$. It is easy to see
also that $(0,1)\setminus \{a \} \subset \Gamma$. The typical local
structure of the level set $\Gamma$ at $t=t_+$ is in Fig.\
\ref{fig:levelsPhi}: there are four arcs of $\Gamma$ emanating from
$t_+$ with equal angles; two of these are $(0,a)$ and $(a, 1)$,
another two should be vertical arcs passing through $t_\pm$ and
ending at $(0,1)$. Since an analogous conclusion is valid in a
neighborhood of $t_-$, we see that necessarily $\Gamma$ contains a
closed Jordan loop passing through $t_\pm$, that we choose as
$\gamma_n $ in \eqref{defI}; see Fig.\
\ref{fig:levelsPhi_dominant_Global}, right.
Analogous arguments allow to determine the orthogonal trajectories
$\Gamma^\perp=\{t \in \C \setminus (-\infty, 0]:\, \Re \Psi
_n(t)=\Re \Psi _n(t_+)\}$, see Fig.\
\ref{fig:levelsPhi_dominant_Global}, left.
Since $\Psi _n(t_-)<\Psi _n(t_+)$, we have that $\Re\Psi _n(t)<\Re
\Psi _n(t_+)$ in the whole connected component of the complement of
$\Gamma^\perp$ containing $t_-$. In particular, $\Re \Psi _n(t)$
attains its strict maximum on $\gamma_n $ at $t_+$, so this is
really the steepest descent curve for $I_n$.

Let us take $\delta >0$ and denote
\begin{equation}\label{def_B_N}
B_n\isdef B_{\delta/\sqrt{n}}(t_+)=\{ t\in \C:\,
|t-t_+|<\delta/\sqrt{n} \}\,, \quad \beta_n=\partial B_n=\{ t\in
\C:\, |t-t_+|=\delta/\sqrt{n} \}.
\end{equation}
We assume $n$ large enough, so that $B_n\subset \D_r\setminus
(-\infty, a]$.  Observe that $t_+$ is the only zero of $\Psi
_n(t)-\Psi _n(t_+)$ in $B_n$, and it is double: in a neighborhood of
$t=t_+$,
\begin{equation*}\label{psiAsymptotics}
\Psi _n(t)-\Psi _n(t_+)=\frac{\Psi ''(t_+)}{2}\, (t-t_+)^2 \left (
1+ \mathcal O (t-t_+)\right)\,, \quad t \to t_+ \,.
\end{equation*}
By \eqref{tau}, $ \Psi _n''(t_+)=2n^{1/2} a^{-3/2}\, \left(1
+o(1)\right) $. Condition $\sqrt{ \Psi _n(t)-\Psi _n(t_+)}>0$ for $t
> t_+$ fixes a single-valued branch of $ \sqrt{
\Psi _n(t)-\Psi _n(t_+)}$ in $B_n$. With this convention, we define
\begin{equation}\label{defZeta}
\zeta_n (t)\isdef   -i n^{1/4} \sqrt{  \Psi _n(t)-\Psi _n(t_+)}\,,
\end{equation}
(where we take the positive square root). It is a conformal mapping
of $B_n$ such that that $\zeta_n\left(\R \cap B_n \right) \subset
i\R$, and $\zeta_n\left(\gamma_n  \cap B_n \right)=(-d_n,
d_n)\subset \R$; moreover, the arc of $\gamma_n  \cap B_n $ in the
upper half plane $\C^+$ is mapped onto the positive semiaxis. Also
\begin{equation*}\label{defZetaBis}
\zeta_n (t)=- i \frac{n^{1/2}}{a^{3/4}}  \, (t-t_+)\left ( 1+
\mathcal O (t-t_+)\right)\,, \quad t \to t_+ \,,
\end{equation*}
and for sufficiently large $n$,
\begin{equation}\label{boundForZeta}
\zeta_n(t)^2 \geq \frac{\delta^2 }{2 a^{3/2}}\,, \quad t\in \gamma_n
\setminus B_n\,.
\end{equation}

\begin{figure}[htb]
\centering
\begin{tabular}{ll}
\hspace{-1.5cm}\mbox{\includegraphics[scale=0.65]{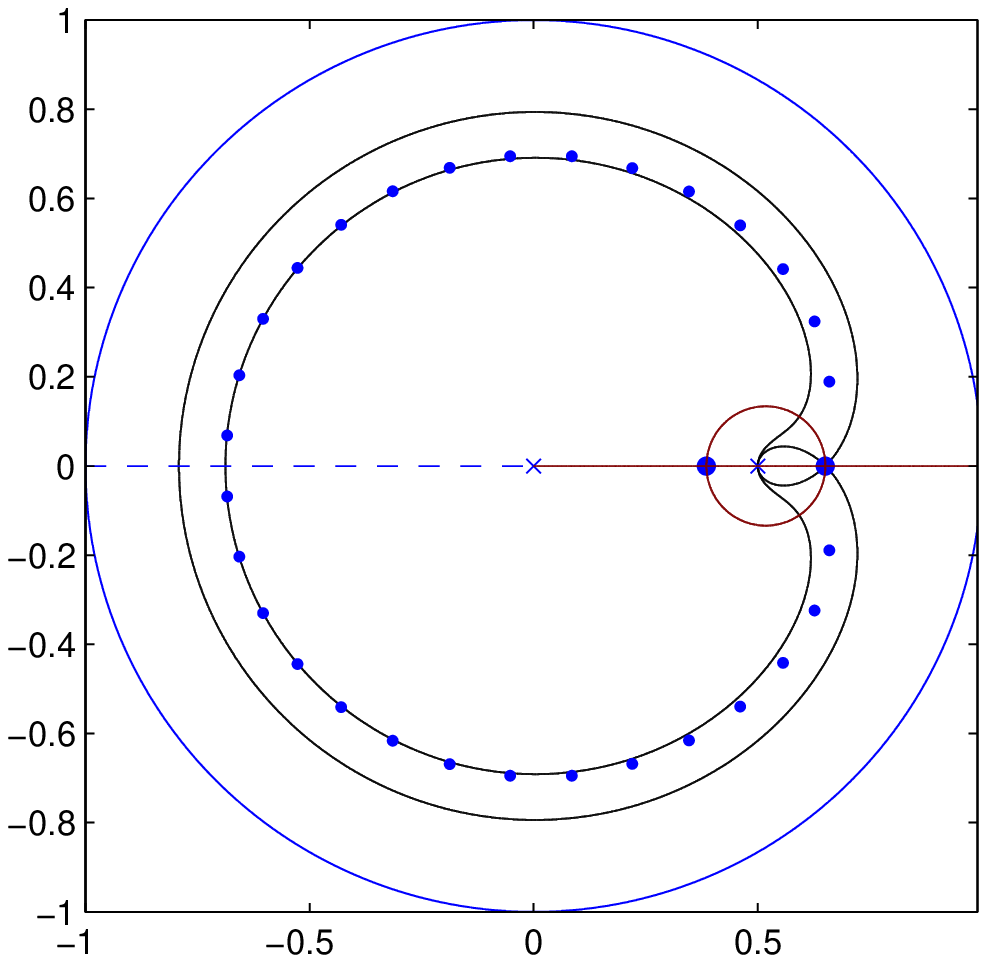}}
&
\hspace{-2cm}\mbox{\includegraphics[scale=0.65]{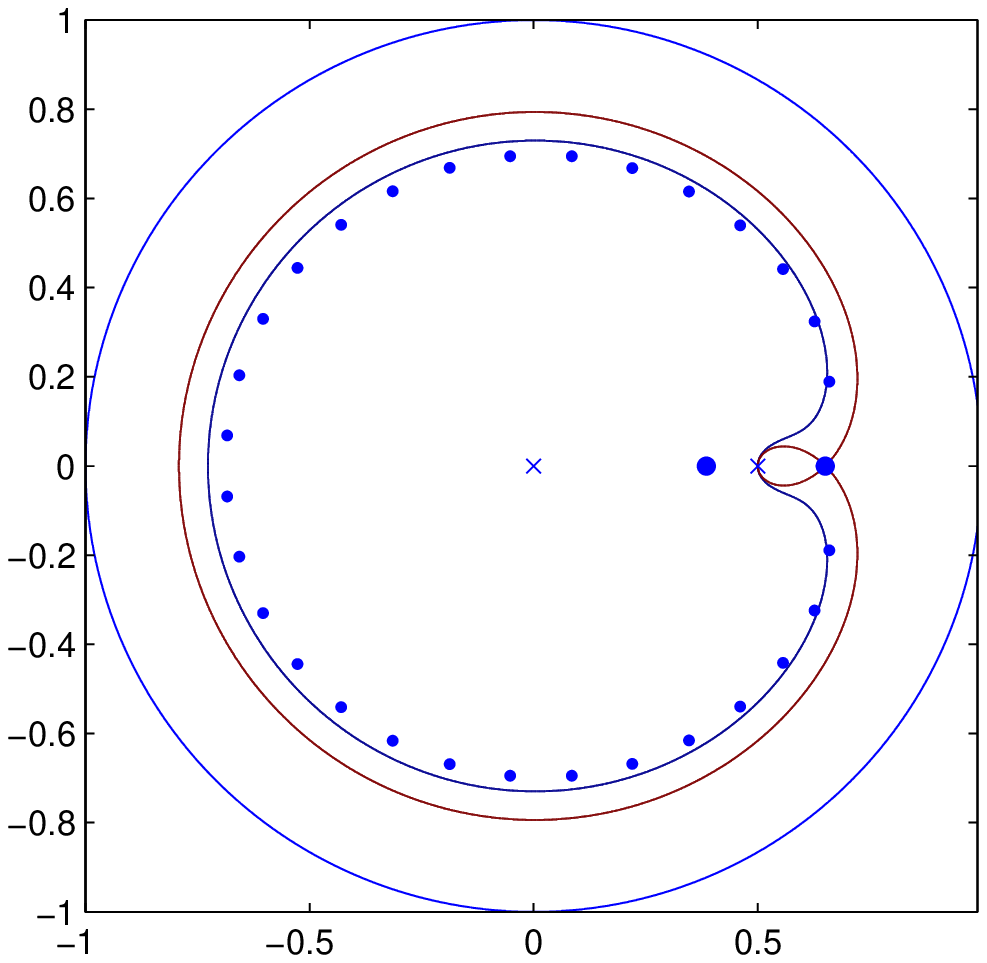}}
\end{tabular}
 \caption{Zeros of $\Phi_{n}$ for $w$ given in \eqref{ess_sing} with
$a=1/2$ and $n=30$, along with the level curves $\Re (\Psi _n(z)-
\Psi _n(t_+))=0$, $\Re ( \Psi _n(z)- \Psi _n(t_+))=\frac{1}{n}\,
\log\left(\frac{1}{2\sqrt{\pi} } \frac{a^{3/4}}{ n^{3/4}} \right)$
and $\Im (\Psi _n(z)- \Psi _n(t_+))=0$ (left). Picture on the right
shows that zeros of $\Phi_{30}$ in a neighborhood of the singular
point $a$ follow more closely the level curve $\Re ( \Psi _n(z)-
\Psi _n(t_+))=\frac{1}{n}\, \log\left(\frac{1}{2\sqrt{\pi} }
\frac{a^{3/4}}{ n^{1/4}} \right)$, as predicted.
}\label{fig:Essential30_dominant}
\end{figure}

For the error function
\begin{equation} \label{defP1}
\begin{split}
P_n(u )  & \isdef - e^{n \Psi _n(t_+)} \, \frac{1}{2\pi i}\,
 \int_{-\infty}^{+\infty} \frac{e^{-\sqrt{n}\, t^2}}{t-u}\, dt \\
 & =  -\frac{1}{2} \, e^{n \Psi _n(t_+) -\sqrt{n} u^2}\,
\mathrm{erfc}\left(- i \, n^{1/4}\, u   \right)\,,
 \quad u\notin \R\,,
 \end{split}
\end{equation}
we define $F_n(z)\isdef P_n(\zeta_n(z))$, $z\in B_n\setminus
\gamma_n $. Then $F_n$ has the following jump across $\gamma_n \cap
B_n$:
$$
F_{n, +}(z)-F_{n, -}(z)=- e^{n \Psi _n(t_+)} \,e^{-\sqrt{n}\,
\zeta_n(z)^2}=-e^{n \Psi _n(z)}\,,
$$
matching the jump of $I_n$ across the same arc. This fact along with
Sokhotsky's formulas show that
\begin{equation}\label{Iexpression}
I_n(z)= \begin{cases} \displaystyle I_n^*(z) + F_n(z)- \frac{1}{2\pi
i}\,
\int_{\beta_n} \frac{F_n(t)}{t-z}\, dt\,, & \text{for } z\in B_n \setminus \gamma_n\,, \\
 \displaystyle  I_n^*(z) - \frac{1}{2\pi i}\, \int_{\beta_n}
\frac{F_n(t)}{t-z}\, dt\,, & \text{for } z\in \C\setminus
(\overline{B_n}\cup \gamma_n^*)\,,
\end{cases}
\end{equation}
where
\begin{equation}\label{E}
I_n^*(z)\isdef -\frac{1}{2\pi i}\, \int_{\gamma_n ^* } \frac{e^{n
\Psi _n(t)}}{t-z}\, dt\,, \quad \gamma_n ^*=\gamma_n \setminus
B_n\,.
\end{equation}
Taking into account \eqref{boundForZeta},
\begin{align*} \label{boundIstar}
 \left|I_n^*(z) \right|& \leq \frac{ e^{n \Psi _n(t_+)}}{\dist (z,\gamma_n ^*)}\, \exp\left(
-\frac{\delta^2 }{2 a^{3/2}}\, n^{1/2}\right)\,, \quad z\notin
\gamma _n^*\,.
\end{align*}
Asymptotics for $F_n$ can be obtained using \cite[Formula
7.1.23]{abramowitz/stegun:1972},
$$
F_n(t)  = \frac{e^{n \Psi _n(t_+) }}{2  \sqrt{\pi}  \,
 } \,    \frac{a^{3/4}}{n^{3/4}(t-t_+)}  \, \left(1+\mathcal O \left(\frac{1}{n^{1/2}
} \right)\right) \,, \quad t\in \beta_n\,;
$$
also for $z\in B_n\setminus \gamma _n$,
$$
\left|\frac{1}{2\pi i}\, \int_{\beta_n} \frac{F_n(t)}{t-z}\,
dt\right|\leq \frac{C n^{-1/4}\, e^{n\Psi_n(t_+)}}{\dist
(z,\beta_n)}\,.
$$
These estimates show that, roughly speaking, in the expressions
\eqref{Iexpression} for $I_n$ function $F_n$ dominates in $B_n$,
while its Cauchy transform along $\beta_n$ does it outside, and in
consequence,
\begin{equation*}\label{I_N_case1}
I_n(z)=\begin{cases}
  \displaystyle  -\frac{t_+^n \FF(t_+) }{t_+-z}  \, \frac{1}{2\sqrt{\pi}
} \left(\frac{a}{ n}\right)^{3/4}\, \left(1 +\mathcal
O\left(\frac{1}{n^{1/2}}\right)\right) + \varepsilon _1(z), &
\text{if } z \in \C\setminus (\overline{B_n}\cup \gamma_n^*), \\
F_n(z) + \varepsilon _2(z), & \text{if } z \in B_n\setminus
\gamma_n\,,
\end{cases}
\end{equation*}
where
$$
|\varepsilon _1(z)|\leq \frac{ t_+^n \FF(t_+)}{\dist (z,\gamma_n
^*)}\, \exp\left( -\frac{\delta^2 }{2 a^{3/2}}\, n^{1/2}\right)\,,
\quad |\varepsilon _2(z)|\leq \frac{C n^{-1/4}\, t_+^n
\FF(t_+)}{\dist (z,C_n)}\,.
$$
This asymptotic expression has several corollaries. First, formulas
\eqref{f_n} show that $f_n^{(1)}$ should have no zeros inside
$\gamma_n $. All zeros lie outside $\gamma_n$ and approach the
circle $\T_a$ with the clock pattern following asymptotically the
solutions of the equation
$$
\left(\frac{z}{t_+}\right)^n =  \frac{1}{2\sqrt{\pi}}\, \left(
\frac{ a}{n}\right)^{3/4} \, \frac{1 }{z-t_+} \,
 \frac{ \FF(t_+)}{\FF(z) }\,.
$$
Away from $t_+$ they should be close to the level curve $\Re ( \Psi
_n(z)- \Psi _n(t_+))=\frac{1}{n}\, \log\left(\frac{1}{2\sqrt{\pi} }
\frac{a^{3/4}}{ n^{3/4}} \right)$ (compare with Fig.\
\ref{fig:Essential30_dominant}, left), although at the boundary
$\beta_n$ of the disk $B_n$ they approach the level curve $\Re (
\Psi _n(z)- \Psi _n(t_+))=\frac{1}{n}\,
\log\left(\frac{1}{2\sqrt{\pi} } \frac{a^{3/4}}{ n^{1/4}} \right)$
(see Fig.\ \ref{fig:Essential30_dominant}, right).

Furthermore, the Verblunsky coefficients $\alpha_n$ for $\Phi_n$
exhibit the following asymptotic behavior:
\begin{equation}\label{verblunskyForEssential}
\alpha_n=-\frac{1}{2\sqrt{\pi}}\, t_+^{n} \FF(t_+)\, \left(
\frac{a}{n}\right)^{3/4} \left( 1+\mathcal O \left(
\frac{1}{n^{1/2}}\right)\right)\,, \quad n \to \infty\,.
\end{equation}
A new feature is the fractional power of $n$ in
\eqref{verblunskyForEssential}, not appearing when the singularities
of $\FF$ on $\T_\rho$ are only poles (cf.\ Corollary
\ref{cor:verblunskyFiniteCase}).

\begin{figure}[htb]
\centering
\begin{tabular}{rl}
\hspace{-1.5cm}\mbox{\begin{overpic}[scale=0.7]{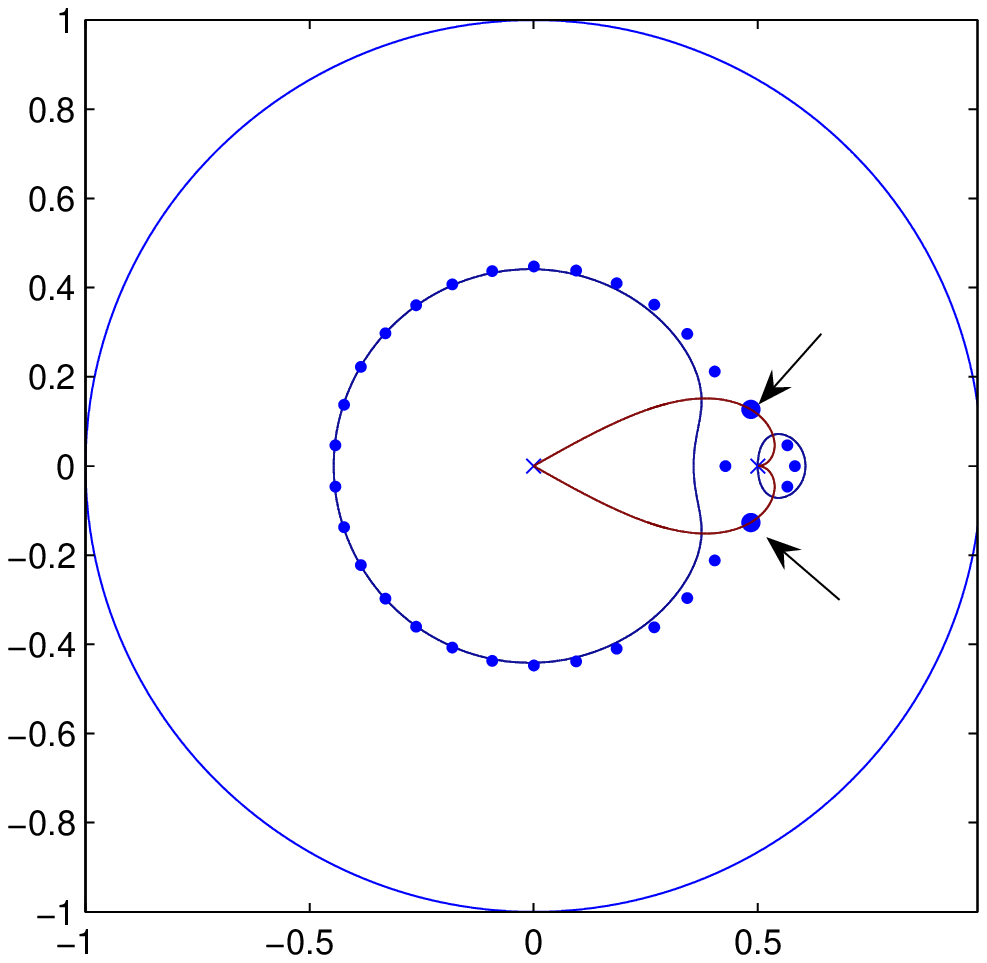}
    \put(72,49){\small $t_+$}
    \put(74,27){\small $t_-$}
   \put(49,37){\small $0$}
\end{overpic}} &
\hspace{-2cm}\mbox{\begin{overpic}[scale=0.7]{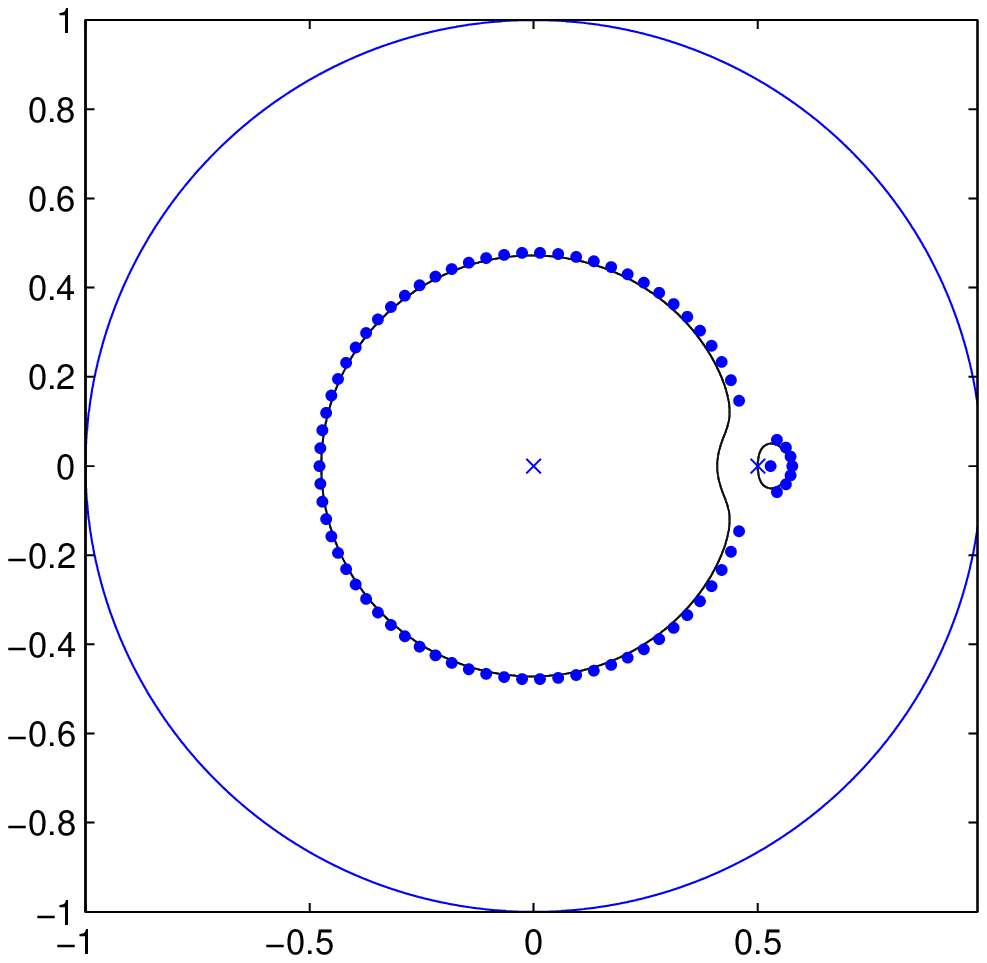}
   \put(49,37){\small $0$}
\end{overpic}}
\end{tabular}
%
%
\caption{Left: zeros of $\Phi_{n}$ for $w$ given in
\eqref{ess_sing2} with $a=1/2$ and $n=30$, contour $\gamma_n$, and
the level curve $\Re ( \Psi _n(z)- \Psi _n(t_+))=\frac{1}{n}\,
\log\left(\frac{1}{2\sqrt{\pi} } \frac{a^{3/4}}{ n^{3/4}} \right)$
(which has two components). Right: same, but with
$n=75$.}\label{fig:Essential}
\end{figure}

Consider now the inverse of the weight function given in
\eqref{ess_sing} by setting
\begin{equation}\label{ess_sing2}
w(t)=\left|\exp\left( \frac{1}{t-a}\right) \right|^2\,, \quad t \in
\T_1\,,
\end{equation}
with $0<a<1$. Surprisingly we observe now a rather different
behavior of the zeros of $\{\Phi_n\}$ (cf.\ Fig.
\ref{fig:Essential}) which is explained by the following analysis.
For $w$ in \eqref{ess_sing2},
\begin{figure}[htb]
\centering
\begin{overpic}[scale=0.8]{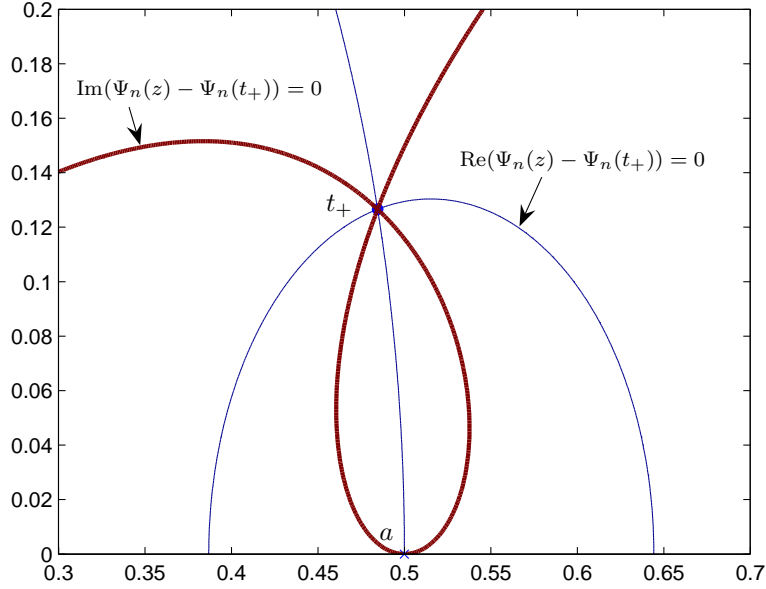}
     \put(15,60){\scriptsize $\Im (\Psi_n(z)-\Psi_n(t_+))=0$}
     \put(58,52){\scriptsize $\Re (\Psi_n(z)-\Psi_n(t_+))=0$}
     \put(43,47){\small $t_+$}
     \put(49,10){\small $a$}
\end{overpic}
 \caption{Level
curves $\Gamma$ (thick lines) and $\Gamma^\perp$ in a neighborhood
of $t_+$ for $a=1/2$ and
$n=30$.}\label{fig:levelsPhi_nondominant_Local}
\end{figure}
$$
\FF(t)=\exp\left( \frac{t}{1- a t }+ \frac{1}{a-t}\right)\,,
$$
and
\begin{align*}
\Psi _n(t) &\isdef  \log t + \frac{1}{n}\, \log \FF(z)=\log t +
\frac{1}{n}\,\left(\frac{t}{1- a t }+ \frac{1}{a-t} \right)
\end{align*}
has saddle points
$$
t_\pm \to a, \quad 1/t_\pm \to 1/a\,,
$$
satisfying
\begin{equation*}\label{tauBis}
  t_+=\overline{t_-}, \quad  t_{\pm}= a \pm i \frac{\sqrt{a}}{\sqrt{n}}+\mathcal O(1/n)\,, \quad n \to
    \infty.
\end{equation*}
Observe that in this case we have no dominant saddle point in a
neighborhood of $a$.

\begin{figure}[htb]
\centering
\begin{tabular}{ll}
\hspace{-1.5cm}\mbox{\begin{overpic}[scale=0.65]{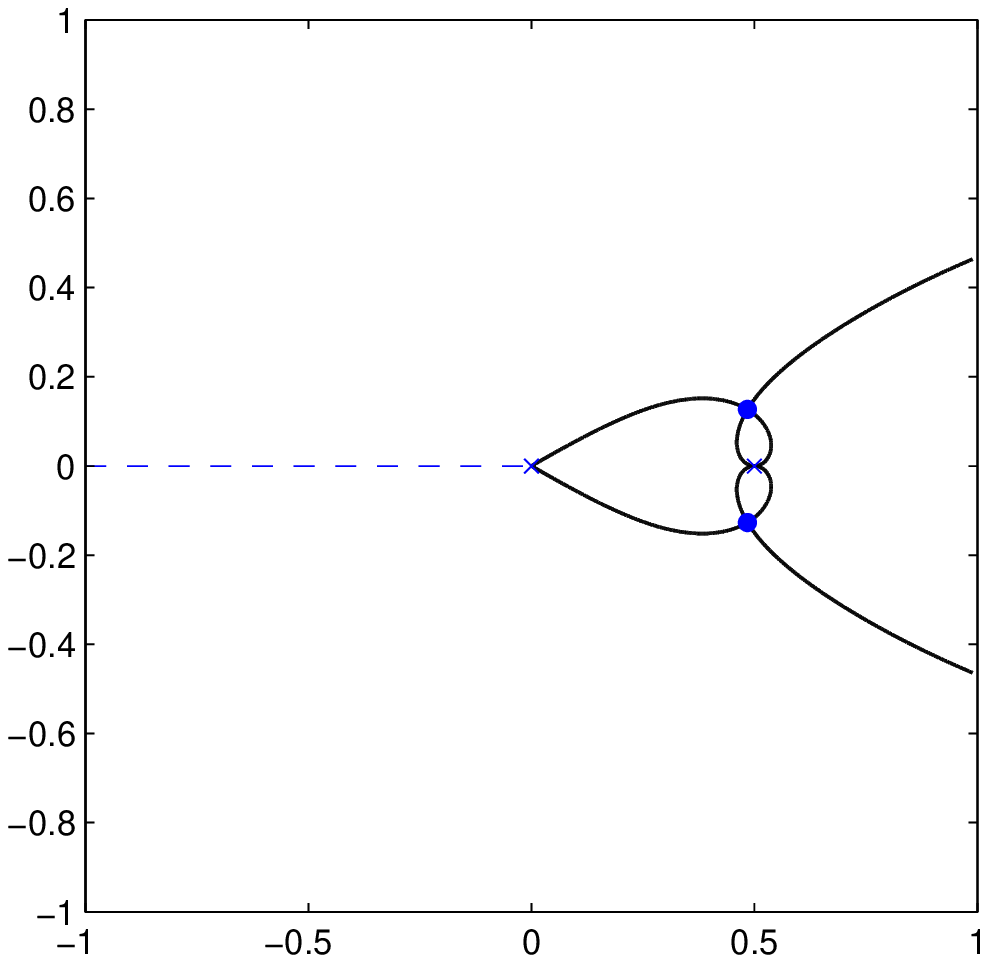}
         \put(51,36){\scriptsize $0$}
         \put(69,38.5){\scriptsize $a$}
\end{overpic}} &
\hspace{-2cm}\mbox{\begin{overpic}[scale=0.65]{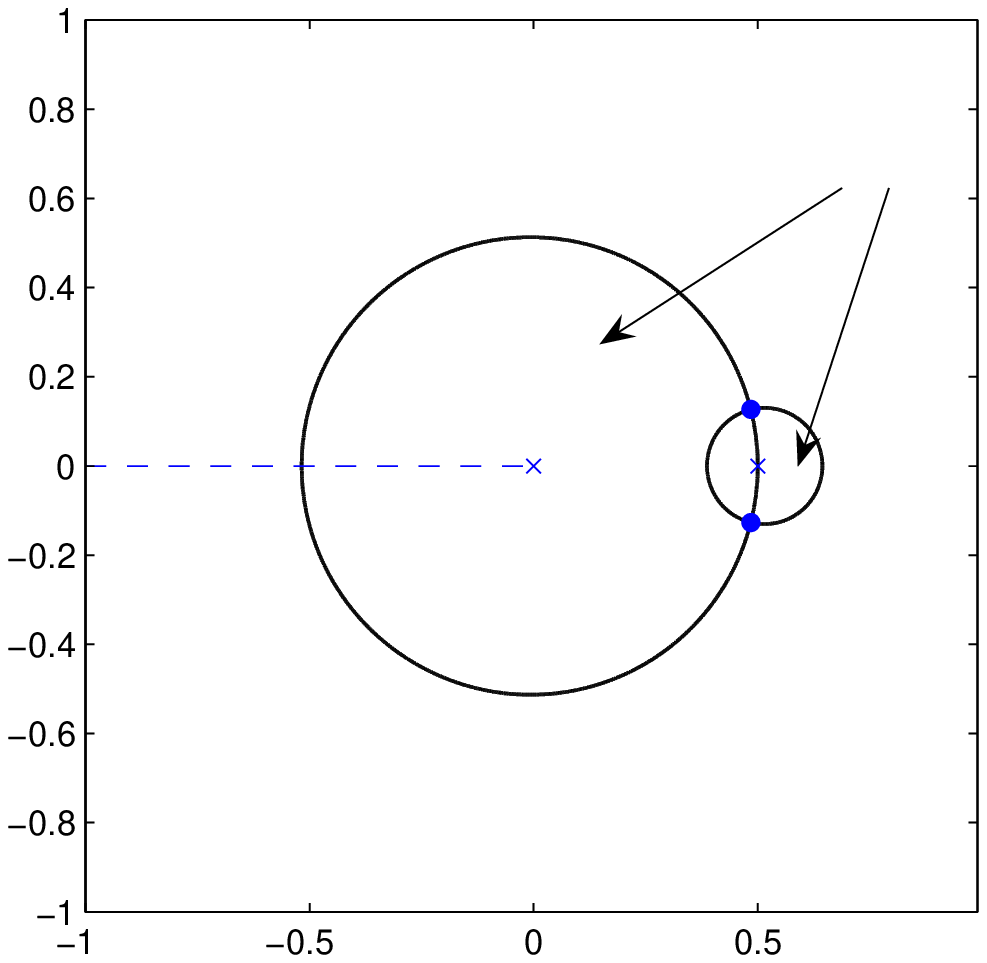}
         \put(51,36){\scriptsize $0$}
         \put(68.5,38.5){\scriptsize $a$}
         \put(67.5,44.5){\scriptsize $t_+$}
         \put(67.5,32.5){\scriptsize $t_-$}
\put(42,59){\small $\Re (\Psi_n(z)-\Psi_n(t_+))<0$}
\end{overpic}}
\end{tabular}
%
%
\caption{Global structure of the level curves $\Gamma$ (left) and
$\Gamma^\perp$ for $a=1/2$ and
$n=30$.}\label{fig:levelsPhi_nondominant_Global}
\end{figure}

An analysis as for weight \eqref{ess_sing} allows us to conclude
that the local structure of the trajectories $\Gamma=\{t\in \C
\setminus (-\infty, 0]:\, \Im (\Psi _n(t)-\Psi _n(t_+))=0\}$ and
$\Gamma^\perp =\{t\in \C \setminus (-\infty, 0]:\, \Re (\Psi
_n(t)-\Psi _n(t_+))=0\}$ is now as depicted in Figure
\ref{fig:levelsPhi_nondominant_Local}, which yields in turn the
global structure is as in Figure
\ref{fig:levelsPhi_nondominant_Global}.
Now we take as $\gamma_n$ in \eqref{defI} the union of arcs of
$\Gamma$ lying in the domain where $\Re (\Psi_n(z)-\Psi_n(t_+))<0$
(see Fig.\ \ref{fig:levelsPhi_nondominant_Global}, right), except
that we slightly depart from $\Gamma$ in an arbitrary small
neighborhood to the right of $a$, where the inequality $\Re
(\Psi_n(z)-\Psi_n(t_+))<0$ is still valid.

Along with the notation introduced in \eqref{def_B_N},
\eqref{defZeta} and \eqref{defP1}, we use and
$B_n^*=\{\overline{z}:\, z\in B_n \}$, $\beta_n^*=\partial B_n^*$,
and
$$
F_n(z)\isdef   P_n(\zeta_n(z)), \; z\in B_n\setminus \gamma_n,
\qquad   F_n^*(z)\isdef P_n(\overline{\zeta_n(\overline{z})}), \; z
\in B_n^*\setminus \gamma_n.
$$
Following the arguments similar to those yielding
\eqref{Iexpression} we obtain
\begin{equation}\label{IexpressionCase2}
I_n(z)= \begin{cases} \displaystyle F_n(z)+I_n^*(z) - \frac{1}{2\pi
i}\, \int_{\beta_n} \frac{F_n(t)}{t-z}\, dt- \frac{1}{2\pi i}\,
\int_{\beta_n^*}
\frac{F_n^*(t)}{t-z}\, dt\,, & \text{for } z\in B_n \setminus \gamma_n\,, \\
\displaystyle F_n^*(z)+I_n^*(z) - \frac{1}{2\pi i}\, \int_{\beta_n}
\frac{F_n(t)}{t-z}\, dt- \frac{1}{2\pi i}\, \int_{\beta_n^*}
\frac{F_n^*(t)}{t-z}\, dt\,, & \text{for } z\in B^*_n \setminus \gamma_n\,, \\
 \displaystyle  I_n^*(z) - \frac{1}{2\pi i}\, \int_{\beta_n}
\frac{F_n(t)}{t-z}\, dt- \frac{1}{2\pi i}\, \int_{\beta_n^*}
\frac{F_n^*(t)}{t-z}\, dt\,, & \text{for } z\in \C\setminus
(\overline{B_n}\cup B^*_n \cup \gamma_n^*)\,,
\end{cases}
\end{equation}
where $I_n^*$ is as in \eqref{E} with $\gamma_n ^*=\gamma_n
\setminus (B_n \cup B_n^*)$.

With contour $\gamma_n$ described above we can perform the
asymptotic analysis of the right hand sides in
\eqref{IexpressionCase2}; again $F_n$ and $F_n^*$ will dominate in
$B_n$ and $B_n^*$, respectively, but the leading term outside these
discs will be given by \emph{the sum} of their Cauchy transforms. In
particular,
$$
I_n(z) \sim - \frac{1}{2\sqrt{\pi} } \left(\frac{a}{ n}\right)^{3/4}
\, \left( \frac{t_+^n \FF(t_+) }{t_+-z}  + \frac{t_-^n \FF(t_-)
}{t_--z}\right)\,, \quad z \in \C\setminus (\overline{B_n}\cup
\overline{B_n^*}\cup \gamma_n^*)\,,
$$
and $I_n(z)$ has (asymptotically) at most one zero inside $\gamma_n
$, originated by the cancelation of these Cauchy transforms. By
\eqref{f_n}, the rest of the zeros of $f_n^{(1)}$ lies outside
$\gamma_n$ and distributes asymptotically as the solutions of the
equation
$$
z^n \FF(z) =   \frac{1}{2\sqrt{\pi} } \left(\frac{a}{
n}\right)^{3/4} \, \left( \frac{t_+^n \FF(t_+) }{t_+-z}  +
\frac{t_-^n \FF(t_-) }{t_--z}\right)\,.
$$
Hence, away from $\gamma_n$ they should be close again to the level
curve $\Re ( \Psi _n(z)- \Psi _n(t_+))=\frac{1}{n}\,
\log\left(\frac{1}{2\sqrt{\pi} } \frac{a^{3/4}}{ n^{3/4}} \right)$
(compare with Fig.\ \ref{fig:Essential}).

Finally, the Verblunsky coefficients $\alpha_n$ of $\Phi_n$ satisfy
$$
\alpha_n=-\frac{1}{\sqrt{\pi}}\, \left( \frac{a}{n}\right)^{3/4}
 \, \Re \left( t_+^{n} \FF(t_+)\right) \, \left( 1+\mathcal O \left(
\frac{1}{n^{1/2}}\right)\right)\,, \quad n \to \infty\,.
$$

\begin{remark}
Although we have considered here the case that $D_{\rm e}$ possesses
a single isolated singularity, it is clear that the same analysis
carries through under the more general assumption that $D_{\rm e}$
can be extended to a domain of the form $\{ |z| > \rho'\}$ with
$\rho' < \rho$ such that $D_{\rm e}$ possesses exactly one isolated
essential singularity of the form \eqref{usePsi} in this domain.
\end{remark}
\begin{remark}
Extensions in various directions of the considerations in this
section are imminently possible; what is really required is $(i)$ a
domain $\{ |z| > \rho'\}$ on which it is known that $D_{\rm e}$
possesses a finite number of singularities and $(ii)$ a description
of the behavior of $D_{\rm e}$ in a vicinity of each singularity.
With this information, a complete asymptotic description of the
function $f_{n}^{(1)}$ is within the reach of careful, and  often
creative, asymptotic analysis of integrals.
\end{remark}

\section*{Acknowledgement} The research of A.M.F.\ was supported, in part, by
a research grant from the Ministry of Science and Technology
(MCYT) of Spain, project code BFM2001-3878-C02, by Junta de
Andaluc\'{\i}a, Grupo de Investigaci\'{o}n FQM229, and by Research Network
``Network on Constructive Complex Approximation (NeCCA)'', INTAS
03-51-6637. A.M.F.\ acknowledges also the support of the Ministry
of Education, Culture and Sports of Spain through the travel grant
PR2003--0104, and the hospitality of the Department of Mathematics
of the Vanderbilt University, where this work was started.

The research of K.T.-R.M.\ was supported, in part, by the National
Science Foundation under grants DMS-0451495 and DMS-0200749.

The research of E.B.S.\ was supported, in part, by the U.S.
National Science Foundation under grant DMS--0296026.

Both A.M.F.\ and K.T.-R.M.\ acknowledge also a partial support of
NATO Collaborative Linkage Grant ``Orthogonal Polynomials: Theory,
Applications and Generalizations'', ref. PST.CLG.979738.

We are grateful to Profs. P.\ Deift, J.\ S.\ Geronimo and P.\ D.\
Miller for very interesting discussions, as well as to Prof.\ B.\
Simon for his tireless activity in this field that stimulated our
research as well.

%

\end{document}